\documentclass[12pt,reqno]{amsart}
\makeatletter
\@namedef{subjclassname@2020}{\textup{2020} Mathematics Subject Classification}
\makeatother
\usepackage{amssymb,amscd,amsthm}
\usepackage[T1]{fontenc}
\usepackage{blkarray}

\usepackage{amsmath,amssymb,graphicx,mathrsfs} 
\usepackage[colorlinks=true,urlcolor=blue,citecolor=red,linkcolor=blue]{hyperref}

\let\frak\mathfrak
\let\Bbb\mathbb

\def\>{\relax\ifmmode\mskip.666667\thinmuskip\relax\else\kern.111111em\fi}
\def\<{\relax\ifmmode\mskip-.333333\thinmuskip\relax\else\kern-.0555556em\fi}
\def\vsk#1>{\vskip#1\baselineskip}
\def\vv#1>{\vadjust{\vsk#1>}\ignorespaces}
\def\vvn#1>{\vadjust{\nobreak\vsk#1>\nobreak}\ignorespaces}

  \let\ssize\scriptstyle
\let\sssize\scriptscriptstyle

\let\Medskip\medskip
\def\medskip{\par\Medskip}
\let\Bigskip\bigskip
\def\bigskip{\par\Bigskip}

\let\Maketitle\maketitle
\def\maketitle{\Maketitle\thispagestyle{empty}\let\maketitle\empty}

\newtheorem{thm}{Theorem}[section]
\newtheorem{cor}[thm]{Corollary}
\newtheorem{lem}[thm]{Lemma}
\newtheorem{prop}[thm]{Proposition}

\numberwithin{equation}{section}

\theoremstyle{definition}
\newtheorem{defn}[thm]{Definition}
\newtheorem{rem}[thm]{Remark}
\newtheorem{prob}[thm]{Problem}
\newtheorem{cex}[thm]{Counterexample}
\newtheorem{ex}[thm]{Example}

\let\mc\mathcal
\let\nc\newcommand

\let\al\alpha
\let\bt\beta
\let\dl\delta
\let\Dl\Delta

\let\gm\gamma
\let\Gm\Gamma

\let\la\lambda
\let\La\Lambda

\let\phi\varphi
\let\si\sigma
\let\Si\Sigma

\let\Ups\Upsilon
\let\om\omega
\let\Om\Omega

\let\der\partial

\let\geq\geqslant

\let\leq\leqslant

\let\on\operatorname
\let\bi\bibitem
\let\bs\boldsymbol

\def\C{{\mathbb C}}
\def\Z{{\mathbb Z}}

\def\Pb{{\mathbb P}}

\def\F{{\mc F}}

\def\+#1{^{\{#1\}}}

\def\diag{\on{diag}}
\def\End{\on{End}}

\def\sln{\mathfrak{sl}_N}

\def\beq{\begin{equation}}
\def\eeq{\end{equation}}
\def\be{\begin{equation*}}
\def\ee{\end{equation*}}

\nc{\bea}{\begin{eqnarray*}}
\nc{\eea}{\end{eqnarray*}}
\nc{\bean}{\begin{eqnarray}}
\nc{\eean}{\end{eqnarray}}
\nc{\bal}{\begin{align*}}
\nc{\eal}{\end{align*}}
\nc{\baln}{\begin{align}}
\nc{\ealn}{\end{align}}

\nc{\Il}{{\mc I_{\bs\la}}}
\nc{\bla}{{\bs\la}}
\nc{\Fla}{\F_\bla}
\nc{\tfl}{{T^*\Fla}}
\nc{\GL}{{GL_n(\C)}}
\nc{\GLC}{{GL_n(\C)\times\C^*}}

\let\sd s 

\def\ddk_#1{\kk_{#1}\<\>\frac\der{\der\<\>\kk_{#1}}}

\def\bul{\mathbin{\raise.2ex\hbox{$\sssize\bullet$}}}
\def\intt{\mathchoice
{\mathop{\raise.2ex\rlap{$\,\,\ssize\backslash$}{\intop}}\nolimits}
{\mathop{\raise.3ex\rlap{$\,\sssize\backslash$}{\intop}}\nolimits}
{\mathop{\raise.1ex\rlap{$\sssize\>\backslash$}{\intop}}\nolimits}
{\mathop{\rlap{$\sssize\<\>\backslash$}{\intop}}\nolimits}}

\let\kk q 
\let\cc c

\let\Ko K

\def\GZ/{Gelfand-Zetlin}
\def\KZ/{{\slshape KZ\/}}
\def\qKZ/{{\slshape qKZ\/}}
\def\XXX/{{\slshape XXX\/}}

\nc{\slnl}{{\sln (\lambda)}}
\nc{\PCN}{{   (\C[x])^N   }}
\nc{\di}{\on{Diag}}
\nc{\dio}{\on{Diag}_0}
\nc{\Mm}{{\mc M}}
\nc{\Nn}{{\mc N}}

\nc{\A}{{\mc C}}

\nc{\PCr}{{  P  (\C[x])^n   }}

\nc{\Pk}{{(\bs{P}^1)^k}}

\nc{\N}{{\Bbb N}}

\nc{\Ll}{{\mc L}}

\nc{\ord}{{\on{ord}\,}}

\nc{\Sing}{{\on{Sing}\,}}
\nc{\sing}{{\on{Sing}\,}}

\nc{\Hess}{{\on{Hess}}}

\nc{\R}{{\Bbb R}}

\let\on\operatorname
\nc{\Kk}{{\bs K}}
\nc{\Ap}{{A_\Phi(z)}}
\nc{\ap}{{A_\Phi(z)}}

\nc{\sv}{{\sing V}}
\nc{\cd}{{\C^n-\Delta}}
\nc{\UT}{{U^0}}   
\nc{\ep}{\epsilon}

\usepackage[all,cmtip]{xy}
\usepackage{empheq}
\usepackage{bm}

\newlanguage\fakelanguage
\newcommand\cyr{\fontencoding{OT2}\fontfamily{wncyr}\selectfont
   \language\fakelanguage}
\DeclareTextFontCommand{\textcyr}{\cyr}

\numberwithin{equation}{section}

\usepackage{calligra}

\DeclareMathOperator{\HOM}{\mathscr{H}\text{\kern -3pt {\calligra\large om}}\,}

\makeatletter
\newsavebox{\@brx}
\newcommand{\llangle}[1][]{\savebox{\@brx}{\(\m@th{#1\langle}\)}%
  \mathopen{\copy\@brx\kern-0.5\wd\@brx\usebox{\@brx}}}
\newcommand{\rrangle}[1][]{\savebox{\@brx}{\(\m@th{#1\rangle}\)}%
  \mathclose{\copy\@brx\kern-0.5\wd\@brx\usebox{\@brx}}}
\makeatother

\newcommand{\bsh}{\begin{shaded}}
\newcommand{\esh}{\end{shaded}}

\newcommand{\dpar}[1]{\mathscr{P}(2,#1)}

\newcommand{\msl}{\{\!\!\{}
\newcommand{\msr}{\}\!\!\}}

\newcommand{\ims}{{\rm Im}}
\newcommand{\Mreg}{\mc M_{\rm reg}}
\newcommand{\Mdiag}{\mc M_{\rm diag}}
\newcommand{\ed}{{\rm d}}
\newcommand{\Igend}{\frak I^{\rm gen}_{\rm d}(\nabla^o)}
\newcommand{\Ifs}{\frak I_{\rm fs}(\nabla^o)}
\newcommand{\Idv}{\frak I_{\rm dv}(\nabla^o)}
\newcommand{\Id}{\frak I_{\rm d}(\nabla^o)}
\newcommand{\dsub}{\mathrel{
\rotatebox[origin=c]{45}{$\subseteq$}}
}
\newcommand{\dsubp}{\mathrel{
\rotatebox[origin=c]{-45}{$\subseteq$}}
}
\newcommand{\DE}[3]{\der_{#1}f_{#2}-\der_{#1}f_{#3}}

\begin{document}
\title[]{On the universality of integrable deformations of solutions of degenerate Riemann--Hilbert--Birkhoff problems}
\author[Giordano Cotti ]{Giordano Cotti$\>^{\circ,\star}$ }
\address{Departamento de Matemática, Faculdade de Ci\^encias da Universidade de Lisboa, Campo Grande Edif\'icio C6, 1749-016 Lisboa, Portugal}
\email{gcotti@fc.ul.pt, gcotti@sissa.it}
\subjclass[2020]{Primary: 34M40, 34M50; Secondary: 32A10. }
\maketitle
\begin{center}
\textit{ $^\circ\>$Departamento de Matemática\\Faculdade de Ci\^encias da Universidade de Lisboa\\
Campo Grande Edif\'icio C6, 1749-016 Lisboa, Portugal\/}
\vskip1,5mm
\textit{ $^\star\>$Grupo de F\'isica Matem\'atica\\Departamento de Matem\'atica, Instituto Superior T\'ecnico\\
Av. Rovisco Pais 1049-001 Lisboa, Portugal\/}
\end{center}
{\let\thefootnote\relax
\footnotetext{\vskip5pt 
\noindent
$^\circ\>$\textit{ E-mail}:  gcotti@fc.ul.pt, gcotti@sissa.it}}

\vskip7mm

\begin{abstract}
This paper addresses the classification problem of integrable deformations of solutions of ``degenerate'' Riemann--Hilbert--Birkhoff (RHB) problems. These consist of those RHB problems whose initial datum has diagonal pole part with coalescing eigenvalues. On the one hand, according to theorems of B.\,Malgrange, M.\,Jimbo, T.\,Miwa, and K.\,Ueno, in the non-degenerate case, there exists a universal integrable deformation inducing (via a unique map) all other deformations \cite{JMU81,Mal83a,Mal83b,Mal86}. On the other hand, in the degenerate case, C.\,Sabbah proved, under sharp conditions, the existence of an integrable deformation of solutions, sharing many properties of the one constructed by Malgrange--Jimbo--Miwa--Ueno \cite{Sab21}. Albeit the integrable deformation constructed by Sabbah is not, {stricto sensu}, universal, we prove that it satisfies a {\it relative} universal property. We show the existence and uniqueness of a maximal class of integrable deformations all induced (via a unique map) by Sabbah's integrable deformation. Furthermore, we show that such a class is large enough to include all generic integrable deformations whose pole and deformation parts are locally holomorphically diagonalizable. In itinere, we also obtain a characterization of holomorphic matrix-valued maps which are locally holomorphically Jordanizable. This extends, to the case of several complex variables, already known results  independently obtained by Ph.G.A.\,Thijsse and W.\,Wasow \cite{Thi85,Was85}.
\end{abstract}

\tableofcontents

\section{Introduction}

\noindent 1.1.\,\,\,{\bf Riemann--Hilbert--Birkhoff (RHB) problems: the analytical and geometrical settings. }
Consider the $n$-dimensional system of ordinary differential equations 
\beq\label{introeq1} z\frac{d}{dz}Y=\mc A(z)Y,\qquad \mc A(z)=z^r\sum_{k=0}^\infty \mc A_kz^{-k},\quad \mc A_0\neq 0,
\eeq
where the series is convergent for $|z|>R_1$, and $r$ is a non-negative integer, called the {\it Poincar\'e rank} of \eqref{introeq1} at $z=\infty$. In 1909, and again in 1913, G.\,Birkhoff addressed the following question \cite{Bir09,Bir13}: 
\vskip2mm
\noindent{\bf Riemann--Hilbert--Birkhoff Problem (analytical version): }{\it Does it exist an analytic matrix-valued function $T(z)=\sum_{k=0}^{\infty}T_kz^{-k}$, with $T_0\in GL(n,\C)$ and the series converging for $|z|>R_2$, such that the transformed equation 
\[z\frac{d}{dz}Z=\widehat{\mc A}(z)Z,\qquad Y(z)=T(z)Z(z),\quad \widehat{\mc A}=T^{-1}\mc AT-zT^{-1}\frac{d}{dz}T,
\]has a polynomial matrix coefficient of the form
\[\widehat{\mc A}(z)=z^r \widehat{\mc A}_0+z^{r-1}\widehat{\mc A}_1+\dots+\widehat{\mc A}_r,\qquad \widehat{\mc A}_j\in M(n,\C)?
\]}
In general, the solvability of the RHB problem is an open problem: although several\footnote{More details on known results will be given in the main body of the paper, see Section \ref{secRHB}.} sufficient conditions for the solvability have been given \cite{Bir13,JLP76,Bal90,Bol94a,Bol94b,BB97}\cite[Ch.\,IV]{Sab07}, it may happen that the RHB Problem does not always admit a positive answer (contrarily to what Birkhoff believed to have proved), see \cite{Gan59,Mas59}.  
\vskip1,5mm
In more geometrical terms, the RHB Problem can be recast as follows:
\vskip1mm
\noindent{\bf Riemann--Hilbert--Birkhoff Problem (geometrical version): }{\it Given a trivial vector bundle $E^{\rm in}$ on a disc $D\subseteq\Pb^1$ centered at $z=\infty$, equipped with a meromorphic connection $\nabla^{\rm in}$ with a pole at $z=\infty$, is it true that $(E^{\rm in},\nabla^{\rm in})$ extends to a pair $(E^o,\nabla^o)$, where $E^o$ is a \emph{trivial} vector bundle on $\Pb^1$, and $\nabla^o$ is a meromorphic connection with only another \emph{logarithmic} pole at $z=0$?}

If one does not insist on the triviality of the vector bundle $E^o$, or on the logarithmic nature of the pole $z=0$ of $\nabla^o$, then the problem is easily solvable, see \cite[App.\,A.2]{Sab98}. What makes the RHB Problem difficult is the conjunction of these requirements.
\vskip2mm
This paper is devoted to the study of {\it families} of RHB problems, rather than focusing on a single one. If $(X,x_o)$ is a pointed complex manifold, we consider families of equations \eqref{introeq1}, parametrized by points of $X$, that is
\beq\label{introeq2}
z\frac{d}{dz}Y=\mc A(z,x)Y,\qquad \mc A(z,x)=z^r\sum_{k=0}^\infty \mc A_k(x)z^{-k},\quad \mc A_k\colon X\to M(n,\C)\text{ holomorphic},
\eeq
with the further assumptions that $\mc A_0$ is not identically zero, and that the series defining $\mc A(z,x)$ is convergent for $|z|>R$ (independent of $x$). The family \eqref{introeq2} can be interpreted as a {\it deformation} of the equation given by the specialization $x=x_o$. 
The deformation is said to be {\it integrable} if there exist holomorphic matrix-valued functions $\Theta_1,\dots,\Theta_{\dim X}$ on $\{|z|>R\}\times X$ such that \eqref{introeq2} is compatible with the system of equations
\beq\label{introeq2.1}
\ed' Y=\Theta Y,\qquad \ed'=\sum_{j=1}^{\dim X}\frac{\der}{\der x^j}\ed x^j,\quad \Theta(z,x)=\sum_{j=1}^{\dim X}\Theta_j(z,x)\ed x^j.
\eeq

In the geometrical description, the joint system of equations \eqref{introeq2},\eqref{introeq2.1} defines a flat meromorphic connection $\nabla$ on a trivial vector bundle $E$ over $D\times X$, with a pole along $\{\infty\}\times X$ of order $r+1$. The connection $\nabla$ is a {\it deformation} of its restriction at $x=x_o$, namely the connection $\iota^*\nabla$ on the pulled-back vector bundle $\iota^*E$, where $\iota\colon D\to D\times X,$ $z\mapsto (z,x_o)$. 
\vskip1,5mm
From now on (including the main part of the paper), we consider\footnote{
We expect that many results mentioned (or obtained) in this paper can be generalized to higher Poincar\'e ranks.} the case $r=1$ only. In this case, the Birkhoff normal form reads
\[\frac{d}{dz}Z=\left(\widehat{\mc A}_0+\frac{1}{z}\widehat{\mc A}_1\right)Z,\qquad \widehat{\mc A}_0,\widehat{\mc A}_1\in M(n,\C).
\]
\vskip2mm
Given a family of RHB problems, parametrized by a pointed complex manifold $(X,x_o)$, we can raise a legitimate question: is solvability an open property? More precisely:
\vskip1,5mm
\noindent{\bf Question: }Assume that the RHB problem defined by the system \eqref{introeq2} 
is solvable when specialized at $x_o\in X$. Is it true that there exists an open neighborhood $U\subseteq X$ of $x_o$ such that the RHB problem is solvable when specialized at $x\in U$?
\vskip1,5mm
According to a theorem originally due to B.\,Malgrange (in \cite[Th.\,2.2.(1)]{Mal83a} the result is formulated with some un-necessary assumptions), and subsequently refined by C.\,Sabbah \cite[Th.\,VI.2.1]{Sab07} (see also \cite[Th.\,5.1(c)]{DH21}),   the answer to the question above is positive if
\begin{itemize}
\item the deformation is integrable, 
\item the manifold $X$ is simply connected\footnote{This can always be realized up to replacing $X$ with a simply connected neighborhood of $x_o$.}.
\end{itemize}
Under these assumptions, one can prove the existence of a sufficiently small open neighborhood $U\subseteq X$ of $x_o$, and the existence and uniqueness of a frame of sections of $E|_{D\times U}$, with respect to which $\nabla|_{D\times U}$ has matrix of connections 1-forms 
\beq\label{introeq3}
\Om(z,x)=-\left(A(x)+\frac{1}{z}B_o\right)\ed z - z\,C(x),\qquad x\in U,
\eeq
where $A$ is a holomorphic matrix-valued function on $U$ (referred to as the {\it pole part} of $\nabla$), $B_o$ is a constant matrix, and $C$ is a holomorphic matrix-valued 1-form on $U$ (referred to as the {\it deformation part} of $\nabla$). See Theorem \ref{thimp}.
\begin{rem}
Assume to be given a family of RHB problems as in \eqref{introeq2} (for the moment not necessarily an integrable deformation), solvable for each value of the deformation parameter $x\in X$. The solutions of the RHB problems define a family of equations 
\beq\label{introeq3.1}
\frac{d}{dz}Z=\left(\widehat{\mc A}_0(x)+\frac{1}{z}\widehat{\mc A}_1(x)\right)Z,\qquad \widehat{\mc A}_0,\widehat{\mc A}_1\colon X\to M(n,\C)\text{ holomorphic.}
\eeq For each fixed $x\in X$, the solutions $Z(z,x)$ are multivalued functions, and they manifest both monodromy and a Stokes phenomenon (at $z=\infty$). The integrability condition of the deformation can be recast in terms of the family of equations \eqref{introeq3.1} only: one can prove that the family of RHB is an integrable deformation if and only if the family of equations \eqref{introeq3.1} is {\it weakly isomonodromic} in the sense of \cite{Guz18}. This means that there exists a fundamental system of solutions $Z(z,x)$ whose monodromy matrix $M:=Z(z,x)^{-1}Z(e^{2\pi\sqrt{-1}}z,x)$ does not depend on $x\in X$. Other refined sets of monodromy data (e.g.\,\,the Stokes matrices) may still depend on $x\in X$.
\end{rem}
\vskip1,5mm
\noindent1.2\,\,\,{\bf Malgrange--Jimbo--Miwa--Ueno universal integrable deformations, and its Sabbah's analogue. } Consider a trivial vector bundle $E^o$ on $\Pb^1$, equipped with a meromorphic connection $\nabla^o$, defined (in a suitable basis of sections) by matrix $\Om_o$ of connection 1-forms in Birkhoff normal form, namely
\beq\label{introeq4}
\Om_o(z)=-\left(A_o+\frac{1}{z}B_o\right)\ed z,\qquad A_o,B_o\in M(n,\C).
\eeq
By the discussion of the previous section, the germ of an arbitrary integrable deformation $(\nabla, E,X,x_o)$ of $(E^o,\nabla^o)$ can be defined by a matrix $\Om(z,x)$ as in \eqref{introeq3}, where $A(x_o)=A_o$.
\vskip2mm
Denote by $\Mreg$ and $\Mdiag$ the subsets of $M(n,\C)$ of regular and diagonalizable $n\times n$-matrices, respectively. These sets can be defined as finite unions of locally closed 
strata in $M(n,\C)$, called {\it bundles of matrices}, introduced by V.I.\,Arnol'd \cite{Arn71}, see Section \ref{secbun}. Each bundle of matrices is defined by fixing the ``type'' of the Jordan form of its elements, i.e. by fixing the number of Jordan blocks of each size (the numerical values of the eigenvalues is free).
\vskip2mm
As it will be explained below, the classification problem of germs of integrable deformations of $\nabla^o$ is of varying difficulty, according to which bundle the matrix $A_o\in M(n,\C)$ belongs (i.e. to the Jordan form of $A_o$).
\vskip2mm
In the case $A_o\in\Mreg$, the classification of germs of integrable deformations is completely understood. According to a theorem of B.\,Malgrange \cite{Mal83a,Mal86}, indeed, there exists a germ of \emph{universal} integrable deformation $(\nabla^{\rm un},\underline{\C}^n,\C^n,\bm u_o)$, where $\underline{\C^n}$ denotes the trivial vector bundle $\C^n\times\C^n$ over $\C^n$. By universality we mean that \emph{any other} germ $(\nabla',E',X',x_o')$ of integrable deformation of $\nabla^o$ is induced by $(\nabla^{\rm un},\underline{\C}^n,\C^n,\bm u_o)$ via a {\it unique} base change, i.e.\,\,via a unique germ of map $\varphi\colon (X',x_o')\to(\C^n,\bm u_o)$. See Theorem \ref{thmal1}.
\vskip2mm
In the case $A_o\in\Mreg\cap\Mdiag$, a more detailed description of the universal integrable deformation $\nabla^{\rm un}$ is available. In this semisimple-regular case, it follows from independent results of B.\,Malgrange and M.\,Jimbo, T.\,Miwa, K.\,Ueno \cite{JMU81,Mal83b,Mal86} that the germ of the universal integrable deformation $\nabla^{\rm un}$ can be defined (in a suitable basis of section of the trivial bundle $\underline{\C^n}$) by the matrix
\begin{multline}\label{introeq5}
\Om_{\rm JMUM}(z,\bm u)=-\left(-\La(\bm u)+\frac{1}{z}\left(\mc B'_o+[\Gm(\bm u),\La(\bm u)]\right)\right)\ed z\\-z\ed\La(\bm u)-[\Gm(\bm u),\ed\La(\bm u)],\qquad \bm u\in\mathbb D
\end{multline}
where
\begin{itemize}
\item $\mathbb D\subseteq \C^n$ is a sufficiently small polydisc centered at $\bm u_o\in\C^n$,
\item $\La(\bm u)={\rm diag}(u^1,\dots, u^n)$, where $\bm u\in\C^n$,
\item $\mc B'_o={\rm diag}(b_1,\dots, b_n)$ is a constant diagonal matrix,
\item $\Gm=(\Gm_{ij})_{ij=1}^n$ is an off-diagonal matrix. 
\end{itemize}
Moreover, these data are such that:
\begin{enumerate}
\item there exists $P\in GL(n,\C)$ such that $P^{-1}A_oP=\La(\bm u_o)$, and the diagonal part of $P^{-1}B_oP$ equals $\mc B'_o$,
\item and $\nabla^{\rm un}$ is formally equivalent, at $z=\infty$, to the matrix connection $-\ed\left(z\La(\bm u)\right)-\mc B'_o\frac{\ed z}{z}$.
\end{enumerate}
\vskip2mm
As long as the matrix $A_o$ is not a regular matrix, the classification problem of germs of integrable deformations of $\nabla^o$ becomes extremely more difficult\footnote{By quoting B.\,Malgrange \cite[Rk.\,3.8]{Mal83a}: \guillemotleft [\dots] le probl\`eme de trouver les solutions [of the integrability equations] passant par un $A_o$ non r\'egulier semble tr\`es difficile.\guillemotright}.
In this paper, we address the classification problem of germs of integrable deformations of $\nabla^o$, in the case $A_o\in\Mdiag$ only, i.e. $A_o$ is a diagonalizable matrices with possibly non-simple spectrum. 
\vskip2mm
In the case $A_o\in\Mdiag$, in \cite{Sab21} C.\,Sabbah considered a connection $\nabla^o$ satisfying the following {\it Property of Partial Non-Resonance} (for short, {\it Property PNR}):
\vskip1,5mm
\noindent{\bf Property PNR: }There exists a matrix $P\in GL(n,\C)$ diagonalizing $A_o$, i.e. $P^{-1}A_oP=\La_o={\rm diag}(u_o^1,\dots, u_o^n)$, and such that the matrix $\mc B_o:=P^{-1}B_oP$ has the following properties\footnote{Here, and in the main body of the text, we denote by $M'$ the diagonal part of a matrix $M$, and by $M''=M-M'$ its off-diagonal part.}:
\begin{itemize}
\item[$(\star)$] $\mc B''_o\in {\rm Im\ ad}(\La(\bm u_o))$.
\item[$(\star\star)$] $\mc B'_o$ is \emph{partially non-resonant}, i.e. we have $(\mc B'_o)_{ii}-(\mc B'_o)_{jj}\notin \Z\setminus\{0\}$ whenever $u_o^i=u_o^j$.
\end{itemize}
\vskip1,5mm
Remarkably, under this assumption, Sabbah proved the existence of a germ of integrable deformation of $\nabla^o$ of the form \eqref{introeq5}, and satisfying the properties (1) and (2) above. For an analytical proof of Sabbah's result, see \cite{Cot21}.
\vskip2mm
Although the integrable deformation constructed by Sabbah shares many properties of the Malgrange--Jimbo--Miwa--Ueno connection $\nabla^{\rm un}$, no claim of universality was formulated in \cite{Sab21}. On the one hand, this is coherent with the fact that, in general, there is no versal deformation of $\nabla^o$ if $A_o$ is not regular (see e.g.\,\,Appendix \ref{appnover}). On the other hand, we will prove that Sabbah's connection still enjoys a \emph{relative} universality property: this was one of the motivating themes of the current work.
\begin{rem}
When $n=2$, the cases $A_o\in\Mreg$ and $A_o\in\Mdiag$ are exhaustive. In such a case a complete classification of germs of integrable deformation has been developed in the very interesting preprint \cite{Her21}. In the current paper, we address the case of arbitrary $n$, and we do not obtain a complete classification as in {\it loc.\,cit.}. Our results, however, still emphasize the richness of the classification in the non-regular case.
\end{rem}
\noindent 1.3.\,\,\,{\bf Results. }The first main result of this paper is of preliminary nature, framed in the general study of operator-valued holomorphic functions. It consists of a characterization of matrix-valued holomorphic functions $A\colon X\to M(n,\C)$, defined on a complex manifold $X$, which are {\it locally holomorphically Jordanizable}. We say that $A$ is locally holomorphically Jordanizable at $x_o\in X$ if there exist an open neighborhood $U\subseteq X$ of $x_o$  and a holomorphic function $P\colon U\to GL(n,\C)$ such that $P^{-1}AP$ is in Jordan canonical form. Clearly, a necessary condition for $A$ to be locally holomorphically Jordanizable at $x_o$ is
\vskip1mm
(I) the existence of a holomorphic function $J\colon U\to M(n,\C)$, in Jordan canonical form for any $x\in U$, and similar to $A(x)$ for any $x\in U$.
\vskip1mm
\noindent The validity of condition (I) only, however, is not sufficient. Set 
\begin{itemize}
\item $\si(A(x))$ to be spectrum (i.e.\,the set of eigenvalues) of $A(x)$ for any $x\in X$;
\item $\la_1,\dots, \la_r\colon U\to \C$ to be the holomorphic eigenvalues functions of $A|_U$ (without counting multiplicities);
\item ${\rm coal}(A)\subseteq X$ to be the {\it coalescence locus} of $A$, namely the set of points $x\in X$ such that 
\[\forall\,V \text{neigh.\,of }x,\,\exists\,z\in V\colon {\rm card}\,\si(A(z))>{\rm card}\,\si(A(x));\]
\item $\mc G_n$ to be the disjoint union of complex Grassmannians of subspaces in $\C^n$, i.e. $\mc G_n=\coprod_{k=0}^nG(k,n)$  (with the complex analytic topology).
\end{itemize}

In Theorem \ref{th1}, we prove that $A$ is locally holomorphically Jordanizable at $x_o$ if and only if conditions (I) above, (II) and (III) below hold true:
\vskip1mm
(II) For each $i=1,\dots, r$, the function $\psi_i\colon U\setminus{\rm coal}(A)\to \mc G_n$ defined by $x\mapsto {\ker}(A(x)-\la_i(x){\rm Id}_n)^n$ admits a limit $L_{i}\in\mc G_n$ at $x=x_o$.
\vskip1mm
 (III) We have $\bigoplus_{i=1}^rL_i=\C^n$.
\vskip1mm
Moreover, standing on results of W.\,Kaballo \cite{Kab76}, we prove that the limit $L_i$, with $i=1,\dots, r$, necessarily equals the space $\mc K[(A-\la_i{\rm Id})^n;x_o]$ of values at $x_o$ of $\C^n$-valued holomorphic functions in the kernel sheaf $\text{\calligra\large ker}\,(A-\la_i{\rm Id})^n$ of the morphism of $\mathscr O_X$-modules
\[(A-\la_i{\rm Id})^n\colon \mathscr O_X^{\oplus n}\to \mathscr O_X^{\oplus n}.
\]We also prove that condition (II) is equivalent to the condition 
\vskip1mm
(II.bis) the function $U\ni x\mapsto \dim \mc K[(A-\la_i{\rm Id})^n;x]$ is continuous at $x_o$, for any $i=1,\dots, r$.
\vskip1mm
Theorem \ref{th1} extends to the case of several complex variables previous results of Ph.G.A. Thijsse and W.\,Wasow, independently obtained in \cite{Thi85,Was85}. See Corollary \ref{corTW}.
\begin{rem}
Let us stress the main differences of our results with those of  \cite{Thi85,Was85}:
\begin{enumerate}
\item The results of both Thijsse and Wasow work under the assumption $\dim_\C X=1$ (more precisely, $X\subseteq \C$ is an open region). In this case, condition (II) above is automatically satisfied. Thijsse realized this standing on the results of \cite{BKL75}. Alternatively, we will deduce this from results of W.\,Kaballo \cite{Kab76,Kab12}. 
\item Wasow's Theorem actually works under the further assumption ${\rm coal}(A)=\emptyset$. In such a case, also condition (III) is automatically satisfied.
\item Thijsse's and Wasow's results are statements of {\it global holomorphic similarity}. If $X\subseteq \C$ is an open region, and if $A$ is locally holomorphically Jordanizable at each point $x\in X$, then one can prove that there exists a globally defined holomorphic matrix $T\colon X\to GL(n,\C)$ such that $T(x)^{-1}A(x)T(x)$ is in Jordan form. This actually is an important peculiarity of all 1-dimensional Stein manifolds $X$, see \cite{Gur88,Lei20,For17} and Theorem \ref{thste1}.
\end{enumerate}
To the best of our knowledge, our characterization of locally holomorphically Jordanizable matrices, depending on several complex variables, was never explicitly formulated in literature. The local holomorphic Jordanizability condition plays a crucial condition in the theory of isomonodromic deformations at irregular singularities with coalescing eigenvalues \cite{CG18,CDG19}. We expect that our characterization will be useful for generalizing results of \cite{CG18,CDG19}, as well as for the study of strata of Dubrovin--Frobenius and flat $F$-manifolds \cite{CG17,CDG20, Cot21b}.
\end{rem}
\vskip2mm
The second main results of the paper concerns the classification problem of integrable deformations of $\nabla^o$ defined by \eqref{introeq4}. We first introduce several classes of germs of integrable deformations of $\nabla^o$, namely:
\begin{itemize}
\item the class $\frak I(\nabla^o)$ of {\it all} germs of integrable deformations of $\nabla^o$;
\item the class $\Id$ of germs of {\it diagonal type} (d-{\it type}) integrable deformations of $\nabla^o$, that is those integrable deformations whose pole and deformation parts ($A$ and $C$ in equation \eqref{introeq3}) are locally holomorphically diagonalizable;
\item the class $\Igend$ of germs of {\it generic} d-type integrable deformations of $\nabla^o$, that is those integrable deformations $(\nabla,E,X,x_o)$ of $\nabla^o$ of d-type whose pole part has holomorphic diagonal form $\Dl_0={\rm diag}(f_1(x),\dots, f_n(x))$ with $\ed_{x_o}f_i\neq \ed_{x_o}f_j$ for $i\neq j$;
\item the class $\Idv$ of germs of {\it diagonal-vainshing type} (dv-{\it type}) integrable deformations of $\nabla^o$, that is those integrable deformations $(\nabla, E,X,x_o)$ which are defined (in a suitable basis) by a matrix of connection 1-forms
\[\widetilde\Om(z,x)=-\left(\Dl_0(x)+\frac{1}{z}\mc B(x)\right)\ed z-z\,\ed \Dl_0(x)+\varpi(x),
\]where the $\Dl_0$ and $\mc B$ are holomorphic matrix-valued functions, and $\varpi$ is a holomorphic matrix-valued 1-form such that
\[\Dl_0(x)={\rm diag}(f_1(x),\dots, f_n(x)),\qquad \mc B''=[\mc L,\Dl_0],\qquad \varpi''=[\ed \Dl_0,\mc L],
\]for a holomorphic matrix-value function $\mc L\colon X\to M(n,\C)$, $\mc L=\mc L''$;
\item the class $\Ifs$ of germs of {\it formally simplifiable} integrable deformations of $\nabla^o$ which are formally equivalent, at $z=\infty$, to the connection 
\[\ed-\ed (z\Dl_0(x))-\frak B'_o\frac{\ed z}{z},
\]where $\Dl_0(x)={\rm diag}(f_1(x),\dots, f_n(x))$ and $\frak B_o'$ is a diagonal constant matrix.
\end{itemize}
In Theorems \ref{mth1} and \ref{mth2} we show that
\begin{enumerate}
\item if $A_o\in\Mreg\cap\Mdiag$, then we have $$\emptyset\neq \Igend\subsetneq \Ifs=\Idv=\Id=\frak I(\nabla^o),$$
\item while if $A_o\in\Mdiag$ and the Property PNR holds true, then we have a richer classification, since
\[\emptyset\neq \Igend\subsetneq \Ifs\subseteq\Idv\subseteq\Id\subseteq\frak I(\nabla^o),
\]where all the inclusions are, in general, expected to be strict (see Remark \ref{rkstrict}). 
\end{enumerate}
For clarity of exposition, let us denote with the unified notation\footnote{Here JMUMS stands for Jimbo--Miwa--Ueno--Malgrange--Sabbah.} $\nabla^{\rm JMUMS}$ the integrable deformation of $\nabla^o$ constructed by Jimbo--Miwa--Ueno--Malgrange (in the case $A_o\in\Mdiag\cap\Mreg$), and by Sabbah (in the case $A_o\in\Mdiag$ \& Property PNR). In Theorem \ref{mth4} we further show the existence and uniqueness of a class $\frak I_{\rm JMUMS}$ of germs of integrable deformations of $\nabla^o$ such that
\begin{enumerate}
\item the integrable deformation $\nabla^{\rm JMUMS}$ is an element of $\frak I_{\rm JMUMS}$,
\item any germ of integrable deformation, which is induced by $\nabla^{\rm JMUMS}$ via a base change map, is an element of $\frak I_{\rm JMUMS}$,
\item $\frak I_{\rm JMUMS}$ is maximal with respect to the properties above.
\end{enumerate}
Moreover, we show that if a germ of integrable deformation is induced by $\nabla^{\rm JMUMS}$ via a base change map, then the germ of such a map is unique. For short, we say that $\nabla^{\rm JMUMS}$ is $\frak I_{\rm JMUMS}$-{\it universal}, see Definition \ref{iun}. 

In the case $A_o\in\Mreg\cap\Mdiag$, we necessarily have $\frak I_{\rm JMUMS}=\frak I(\nabla^o)$ by universality of the JMUM integrable deformation. In the more general case $A_o\in\Mdiag$ \& Property PNR, we show that $\frak I_{\rm JMUMS}$ is large enough to satisfy the inequality
\beq\label{introeq6}
\Igend\subseteq \frak I_{\rm JMUMS}\subseteq \Ifs.
\eeq 
In addition to that, in Appendix \ref{nPNR} we show that the validity of the PNR condition is a  sharp condition for the results above. More precisely, we exhibit a connection $\nabla^o$ with $A_o\in\Mdiag$ and for which the Property PNR cannot hold true: we show that the germs in $\Igend$ cannot be induced by a single germ of integrable deformation, in opposition to \eqref{introeq6}.
\vskip2mm
The proof of the left inequality in \eqref{introeq6} is based on a remarkable initial value property of an overdetermined system of non-linear PDEs, that we call {\it generalized Darboux--Egoroff equations}. Consider $n$ holomorphic functions $f_1(\bm x),\dots, f_n(\bm x)$ in $d$ complex variables $\bm x=(x^1,\dots, x^d)$. Let $b_1,\dots, b_n\in\C$ be arbitrary constants. The generalized Darboux--Egoroff system ${\rm DE}_{d,n}((f_i)_{i=1}^n;(b_i)_{i=1}^n)$, in the $n^2-n$ unknown functions $(F_{kh}(\bm x))_{k,h=1}^n$, with $k\neq h$, is given by
\begin{multline*}
(\DE{j}{h}{k})\der_iF_{kh}-(\DE{i}{h}{k})\der_jF_{kh}=\\
\sum_{\ell=1}^n(\DE{i}{\ell}{k})(\DE{j}{h}{\ell})F_{k\ell}F_{\ell h}-\sum_{\ell=1}^n(\DE{j}{\ell}{k})(\DE{i}{h}{\ell})F_{k\ell}F_{\ell h},
\end{multline*}
\begin{multline*}
(f_h-f_k)\der_iF_{kh}=(b_h-b_k-1)(\DE{i}{h}{k})F_{kh}\\+\sum_{\ell=1}^n(\DE{i}{\ell}{k})(f_h-f_\ell)F_{k\ell}F_{\ell h}
-\sum_{\ell=1}^n(f_\ell-f_k)(\DE{i}{h}{\ell})F_{k\ell}F_{\ell h},
\end{multline*}
for any $i,j=1,\dots, d$, and any $k,h=1,\dots, n$, with $k\neq h$. The solutions $F_{kh}$ of this system can be arranged in an off-diagonal matrix $F\colon U\subseteq\C^d\to M(n,\C)$.

Theorem \ref{mth3}, which is the third main result of the paper, asserts that
\begin{itemize}
\item if $\bm x_o\in\C^d$,
\item if $\ed_{\bm x_o}f_h\neq \ed_{\bm x_o} f_k$ for any $h,k=1,\dots, n$, with $k\neq h$,
\item if $b_h-b_k\not\in\Z\setminus\{0\}$ whenever $f_h(\bm x_o)=f_k(\bm x_o)$,
\end{itemize}
then any formal power series solution $F$ of ${\rm DE}_{d,n}((f_i)_{i=1}^n;(b_i)_{i=1}^n)$, centered at $\bm x_o$,  is uniquely determined by its initial value $F(\bm x_o)$. Consequently, the same result holds for any analytic solution defined in a neighborhood of $\bm x_o$. 
This generalizes the results of \cite[Lemma 5.21]{Cot21} \cite[Lemma 6.16]{Cot21b} for the (standard) Darboux--Egoroff equations.
\vskip1,5mm
\noindent 1.4\,\,\,{\bf Structure of the paper. } In Section \ref{secbun}, we review the theory of bundles of matrices, introduced by V.I.\,Arnol'd. Special emphasis is given on both geometrical and combinatorial aspects of the theory. In particular, we review its connection with the theory of double partitions of integers (an aspect which is missing in the original source \cite{Arn71}), as well as some of the main results on the decomposition 
of the space of matrices defined by the bundles.

Section \ref{secJ} is devoted to the study of matrix-valued holomorphic maps which are locally/globally holomorphically similar. After reviewing known results due to R.M.\,Guralnick, J.\,Leiterer, and F.\,Forstneri\v{c}, we address the study of Jordan forms of matrix-valued holomorphic maps. The first main result of the paper, Theorem \ref{th1}, is formulated and proved. 

In Section \ref{secRHB}, we first review known results on families of RHB problems due to B.\,Malgrange, M.\,Jimbo, T.\,Miwa, K.\,Ueno, and C.\,Sabbah. We subsequently introduce several classes (d/dv/fs-types) of germs of integrable deformations of a solution of a RHB problem, and we study their inclusive relations, see Theorems \ref{mth1} and \ref{mth2}. We introduce the notion of $\frak I$-universality, and we prove that the Sabbah's integrable deformation satisfies a relative universal property for a suitable maximal class $\frak I_{\rm JMUMS}$, see Theorem \ref{mth4}. Furthermore, we study the generalized Darboux--Egoroff system of PDEs, and we prove its initial value property, Theorem \ref{mth3}.

In Appendix \ref{appmal}, we recall Malgrange's proof of existence of a universal integrable deformation for the connection \eqref{introeq4} with $A_o\in\Mreg$. 

In Appendix \ref{appnover}, we show via an example that if $A_o\not\in\Mreg$, then in general there is no versal integrable deformation of \eqref{introeq4}. 

In Appendix \ref{nPNR}, we show via an example that the Property PNR is a sharp condition for the existence of an integrable deformation inducing all germs of integrable deformations of generic d-type.
\vskip2mm
\noindent {\bf Acknowledgements. }The author is thankful to R.\,Conti, G.\,Degano, D.\,Guzzetti, C. Hertling, 
P.\,Lorenzoni, 
D.\,Masoero, A.T.\,Ricolfi, and 
C.\,Sabbah 
for several valuable discussions. This research was supported by the FCT Projects PTDC/MAT-PUR/ 30234/2017 ``Irregular connections on algebraic curves and Quantum Field Theory'', UIDB/00208/2020 (DOI: 10.54499/UIDB/00208/2020), UIDP/00208/2020 (DOI: 10.54499/UIDP/00208/2020), 2021.01521.CEECIND (DOI: 10.13039/501100005765), and the project 2022.03702.PTDC  (DOI: 10.54499/2022.03702.ptdc) GENIDE - ``GEometry of Non-generic Isomonodromic DEformations''. The author is a member of the COST Action CA21109 CaLISTA.

\section{Bundles of matrices}\label{secbun}

In what follows by {\it multiset} we mean a pair $(A,m)$ where $A$ is a finite set and $m\colon A\to \mathbb N^*$ is an arbitrary function, called {\it multiplicity function}. More informally, a multiset consists of a finite collection of objects (the elements) which may occur more than once: the element $a\in A$ will occur exactly $m(a)$ times.

For short, we will represent multisets by listing their elements, with multiplicity, between $\msl\dots \msr $ brackets.
\vskip2mm
Partitions of integers provide examples for multisets.
Given a non-negative number $n$, a multiset $\la=\msl \la_1,\dots, \la_r\msr$ of non-negative integers is a partition of $n$ provided that $n=\sum_{i=1}^r\la_i$. We denote by $\mathscr P(n)$ the set of partitions of $n$.

\subsection{Double partitions} A {\it double partition} of a positive integer $n$ is a double sum representation of $n$ as follows 
\begin{align}
n&=n_1+\dots+n_k,\qquad n_j=b_{j1}+\dots+b_{jl_j},\quad j=1,\dots,k,
\end{align}
with $b_{ij}\in\Z_{>0}$, and where $k,l_1,\dots,l_k$ are arbitrary positive integers. The numbers $n_1,\dots,n_k$ are the {\it rough parts} of the double partition, whereas the numbers $b_{ij}$ are its {\it fine parts}. The order of rough parts, and of fine parts --inside a single rough part-- is not relevant. 
\vskip2mm
A double partition of $n$ can thus be identified with the datum of a multiset of ordinary partitions 
of the summands $n_i$'s, with $i=1,\dots, k$. We denote by $\dpar{n}$ the set of double partitions of $n$, and by $p(2,n)$ its cardinality.
\vskip2mm
For short, we use the notation $\bm b=\msl b_{11},\dots, b_{1l_1};\dots;b_{k1},\dots, b_{kl_k} \msr$ for the multiset of partitions. The integer $k$ is the {\it rough length} of $\bm b$, denoted by $|\!|\bm b|\!|$, whereas the integers $l_1,\dots,l_k$ are the {\it fine lengths} of $\bm b$.
\vskip2mm
If $\bm b=\msl b_{11},\dots, b_{1l_1};\dots;b_{k1},\dots,b_{kl_k}\msr \in\mathscr P(2,n),$ without loss of generality we may (and will) assume that $\bm b$ is {\it monotonically ordered}: this means that $b_{j1}\geq b_{j2}\geq \dots\geq b_{jl_j}$ for any $j=1,\dots, k$.

\subsection{Bundles of matrices} Given $n\in\N^*$, we denote by $M(n,\C)$ the complex vector space of $n\times n$-matrices. 
\vskip2mm
Two matrices $A_1,A_2\in M(n,\C)$ have the {\it same Jordan type} if their Jordan forms differ by their eigenvalues only, the number of distinct eigenvalues and the orders of the Jordan blocks being the same. 
A {\it bundle of matrices} is a maximal set of matrices with the same Jordan type.
\vskip2mm
Bundles of matrices are in 1-1 correspondence with double partitions of $n$.  
\vskip2mm
Given the double partition $\bm \la=\msl \la_{11},\dots, \la_{1l_1};\dots;\la_{k1},\dots,\la_{kl_k}\msr$, we define $\mc M_{\bm \la}$ to be set the of matrices with Jordan form 
\beq\label{JF}
\bigoplus_{i=1}^k\bigoplus_{j=1}^{l_i} J_{\la_{ij}}(\mu_i), \qquad J_h(\mu)=\begin{pmatrix}
\mu&1&0&\dots\\
0&\mu&1&\dots\\
\vdots&&\ddots&\\
0&\dots&&\mu\\
\end{pmatrix}\in M(h,\C).
\eeq
The partition $\msl \la_{i1},\dots,\la_{il_i}\msr$ is called the {\it Segre characteristic} of the eigenvalue $\mu_i$.
\vskip2mm
The space $M(n,\C)$ admits thus a decomposition into bundles,
\[
M(n,\C)=\coprod_{\bm \la\in\mathscr P(2,n)}{{\mc M}}_{\bm \la}. 
\]

We will also adopt the following 
notation: the bundle $\mc M_{\bm \la}$, associated with ${\bm \la}\in\dpar{n}$ will be labelled by the string
\[\mu_1^{\la_{11}}\mu_1^{\la_{12}}\dots \mu_2^{\la_{21}}\mu_2^{\la_{22}}\dots\mu_k^{\la_{k1}}\mu_k^{\la_{k2}}\dots
\]where $\mu_1,\dots,\mu_k$ denote distinct eigenvalues.
\vskip2mm
For example, for $n=2$, we have $3$ bundles of matrices, labelled by
\[\mu_1\mu_2,\quad \mu_1^2,\quad \mu_1\mu_1,
\]which correspond to the Jordan forms
\[\begin{pmatrix}
\mu_1&0\\
0&\mu_2
\end{pmatrix},\qquad 
\begin{pmatrix}
\mu_1&1\\
0&\mu_1
\end{pmatrix},\qquad
\begin{pmatrix}
\mu_1&0\\
0&\mu_1
\end{pmatrix}
\]respectively.
\vskip2mm
For $n=3$, we have $6$ bundles of matrices, labelled by
\[
\mu_1\mu_2\mu_3,\quad  \mu_1^2\mu_2,\quad \mu_1\mu_1\mu_2,\quad \mu_1^3,\quad \mu_1^2\mu_1,\quad \mu_1\mu_1\mu_1.
\]
\subsection{Number of bundles} 
By the preceding paragraph, the number of bundles of matrices in $M(n,\C)$ equals the number $p(2,n)$ of double partitions of $n$. The first values of $p(2,n)$, for $n=1,\dots, 20$, are 
\begin{multline*}
1,\,3,\,6,\,14,\,27,\,58,\,111,\,223,\,424,\,817,\,1527,\,2870,\,5279,\,9710,\,17622,\\
\,31877,\,57100,\,101887,\,180406,\,318106,\dots
\end{multline*}

The history of this numerical sequence is quite rich and interesting. In 1854, A.\,Cayley first introduced and studied the sequence of numbers $p(2,n)$, see \cite{Cay55}. He claimed that the same numbers arise in the classification of homographies of the projective space $\Pb^n$. Furthermore, Cayley noticed that the first numbers of this sequence (minus 2, and with some computational mistakes for $n=5,7,8$) appeared in \cite{Syl51}, a study of intersections of quadrics by J.J.\,Sylvester. 
Subsequently, all these research directions were extensively developed by C.\,Segre. In 1883, C.\,Segre completed the classification of intersection of two quadrics in a projective space, one of the main topics of his {\it Tesi di Laurea} \cite[Pt.\,II, \S3]{Seg83}. Moreover, in \cite[Ch.\,II.14]{Seg12} Segre presented a complete classification of collineations in projective spaces. In both classification problems, a multiset of partitions (the {\it Segre symbol}) is a classifying invariant. Segre symbols are indeed double partitions of integers, and their number thus equal $p(2,n)$. 
For further details see \cite[Book II, Ch.\,VIII]{HP94a} \cite[Book IV, Ch.\,XIII, \S10--11]{HP94b}. See also \cite{Bel16, FMS21}.
\vskip2mm
Introduce the generating function 
\[\mc P(2;z):=\sum_{n=0}^\infty p(2,n)z^n,\quad p(2,0):=1.
\]Denote by $p(n)$ the number of ordinary partitions of $n$. 
\begin{thm}[\cite{Cay55}]\label{thmgenf}
We have $$\mc P(2;z)=\prod_{n=1}^\infty(1-z^n)^{-p(n)}.$$
\end{thm}
\proof For $n\in\N$, the number $p(2,n)$ equals the number of representations of $n$ as following sums
\[n=\sum_{k=1}^\infty k(m_{k,1}+\dots + m_{k,p(k)}),\quad m_{i,j}\geq 0.
\]This follows from the interpretation of a double partition as a multiset of partitions, and the $m_{i,j}$'s represent the multiplicities. 
We have
\[
\pushQED{\qed} 
\prod_{n=1}^\infty(1-z^n)^{-p(n)}=\prod_{n=1}^\infty\sum_{\substack{m_{j}=0\\ 1\leq j\leq p(n)}}^\infty z^{n(m_{1}+\dots+m_{p(n)})}=\sum_{\substack{m_{i,j}=0\\ 1\leq i,j}}^\infty z^{1m_{1,1}+2(m_{2,1}+m_{2,2})+\dots}.\qedhere
\popQED
\]
\begin{cor}
We have the recursive formula
\[p(2,n)=\frac{1}{n}\sum_{k=1}^n\si(k)p(2,n-k),\quad \si(k):=\sum_{d|k}d\cdot p(d).
\]
\end{cor}
\proof
From Theorem \ref{thmgenf}, we obtain $\log \mc P(2,z)=\sum_{n,m=1}^\infty \frac{p(n)}{m}z^{nm}$. The recursive formula follows from the identity $\frac{d}{dz}\mc P(2,z)=\mc P(2,z)\frac{d}{dz}(\log \mc P(2,z))$.
\endproof

The numbers $p(2,n)$ have an exponential growth, as described by the following result, due to R.\,Kaneiwa and V.\,M.\,Petrogradsky, which provides an analog of the well-known Hardy--Ramanujan asymptotic formula for $p(n)$, see \cite{HR18}.
\begin{thm}[{\cite{Kan79,Kan80,Pet99}}]
We have the following asymptotic expansion
\[
\pushQED{\qed}
\ln p(2,n)=\left(\frac{\pi^2}{6}+o(1)\right)\frac{n}{\ln n},\quad n\to\infty.\qedhere
\popQED
\]
\end{thm}
\begin{rem}
Although \cite{Kan79,Kan80} provide more precise terms of the expansion, the results of \cite{Pet99} have wider implications. Ordinary partitions and double partitions are just the first instances of $r$-{\it fold partitions}, defined as representations of an integer $n\in\N$ as $r$-fold sums of non-negative integers, see \cite{Kan79}. Denoting by $p(r,n)$ the number of $r$-fold partitions of $n$, with $p(r,0):=1$ for any $r\geq 1$, the generating function $\mc P(r,z):=\sum_np(r,n)z^n$ satisfies the identity 
\[\mc P(r,z)=\prod_{n=1}^\infty(1-z^n)^{-p(r-1,n)},\quad r\geq 2.
\]See \cite{Kan79}. As an application of \cite[Th.\,2.1]{Pet99}, one obtains that
\[\ln p(r,n)=\left(\frac{\pi^2}{6}+o(1)\right)\frac{n}{\ln^{(r-1)} n},\quad n\to\infty,
\]where $\ln^{(k)}x:=\underbrace{\ln\ln\dots\ln}_{k\text{ times}} x$ for $k\geq 1$. More general applications of \cite[Th.\,2.1]{Pet99} allows an estimate of the growth of the number of some generalized partitions, and growth of free polynilpotent finitely generated Lie algebras. See also \cite{Pet00}.
\end{rem}

\subsection{Bundles as fibered spaces}  

Given $\bm a=(a_1,\dots,a_h)\in\N^h$, with $h\geq 1$, set $|\bm a|:=\sum_{j=1}^ha_j$. Let $\frak S_{\bm a}$ to be the subgroup of the symmetric group $\frak S_{|\bm a|}$ defined by $\frak S_{\bm a}:=\frak S_{a_1}\times\dots\times\frak S_{a_h}$, where $\frak S_{a_i}$ is the symmetric group on the elements $\{(\sum_{k=1}^{i-1}a_k)+1,(\sum_{k=1}^{i-1}a_k)+2,\dots, \sum_{k=1}^{i}a_k\}$, for $i=1,\dots, h$.
\vskip2mm
Define the {\it configuration space $\mc C_{\bm a}$ of $|\bm a|$ colored points in the plane} as the quotient
\[\mc C_{\bm a}:=\left(\C^{|\bm a|}\setminus \Dl\right)/\frak S_{\bm a},
\]where $\Dl$ be the union of {big diagonal} hyperplanes in $\C^n$, defined by the equations
\[\Dl:=\bigcup_{i<j}\{\bm u\in\C^n\colon u^i=u^j\}.
\]The tuple $\bm a$ dictates the coloring of the points: in total we have $h$ colors, and for each $i=1,\dots,h$ we have $a_i$ points with the same $i$-th color. The order of points with the same color is not relevant.
\vskip2mm
The bundles $\mc M_{\bm \la}$ are fibered spaces over suitable colored configuration spaces $\mc C_{\bm a}$. More precisely, for each $\bm \la\in\mathscr P(2,n)$ denote by $\bm m_{\bm \la}=(m_1,\dots,m_h)\in \N^h$ the tuple of multiplicities of the elements of the multiset $\bm \la$. We have a natural map $\pi_{\bm \la}\colon\mc M_{\bm \la}\to \mc C_{\bm m_{\bm \la}}$, defined by associating to a matrix $A$ its {\it colored spectrum}, i.e.\,\,the set of its {\it colored distinct eigenvalues}. Two distinct eigenvalues have the same color if they have the same Segre characteristics. The fibers of $\pi_{\bm \la}$ are similarity orbits of the Jordan matrices \eqref{JF}.

\begin{thm}[\cite{Arn71}]
The codimension of the bundle $\mc M_{\bm \la}$ equals
\[c=\left[\sum_{j=1}^{|\!|\bm \la|\!|}\la_{j1}+3\la_{j2}+5\la_{j3}+\dots\right]-|\!|\bm \la|\!|.
\]
\end{thm}
\proof
The dimension of $\mc M_{\bm \la}$ equals the dimension of the base space $\mc C_{\bm m_{\bm \la}}$ plus the dimension of the fiber. We have $\dim \mc C_{\bm m_{\bm \la}}=|\!|\bm \la|\!|$, and the codimension of the similarity orbit of the Jordan matrix \eqref{JF} equals the first summand above. See \cite[Ch.\,VIII]{Gan59}.
\endproof

\subsection{Stratification of bundles} Consider the decomposition of the space of matrices $M(n,\C)$ into bundles, that is $M(n,\C)=\coprod_{\bm \la\in\mathscr P(2,n)}\mc M_{\bm \la}$. The following theorem summarizes some of the main properties of this decomposition. For detailed proofs see \cite[Sec.\,4.1]{DTD22}.
\begin{thm}[\cite{DTD22}]\label{dtd22}The following properties hold.
\begin{enumerate}
\item Each bundle $\mc M_{\bm \la}$ is locally closed in the standard topology of $M(n,\C)$.
\item For each double partition $\bm\la$, there exists a finite number of distinct double partitions $\bm\la_1,\dots,\bm\la_k$, with $\bm\la_1=\bm\la$, such that
\[\overline{\mc M_{\bm \la}}=\bigcup_{j=1}^k\mc M_{\bm\la_j},\quad\text{and}\quad\overline{\mc M_{\bm\la_j}}\subsetneq \overline{\mc M_{\bm\la}},\quad j=2,3,\dots,k.
\]In particular, each closure $\overline{\mc M_{\bm \la}}$ is a constructible set of $M(n,\C)$.
\item If $\mc M_{\bm \la}\cap \overline{\mc M_{\bm \nu}}\neq\emptyset$ then $\mc M_{\bm \la}\subseteq \overline{\mc M_{\bm \nu}}$.\qed
\end{enumerate}
\end{thm}
\begin{rem}
In the original source \cite[\S 5.3]{Arn71}, the decomposition of the space of matrices into bundles is said to define a finite {\it semi-algebraic stratification}: each bundle is claimed to be a \guillemotleft\,semi-algebraic smooth submanifold\guillemotright\, of $M(n,\C)$. To the best of our knowledge, despite the efforts of several authors studying the topology of the bundles decompositions, no detailed arguments for these statements have been provided. See for example \cite[pag.\,670]{EEK99}, where the authors claim the space of matrices to be a {\it stratified manifold} in the sense of \cite[pag.\,41]{AGV85}, by implicitly assuming the validity of point (2) of Theorem \ref{dtd22} above. Also a complete proof of the fact that each bundle is open in its closure (i.e.\,locally closed) was missing. The recent preprint \cite{DTD22} fills this gap in the literature. Notice that point (1) of Theorem \ref{dtd22} works in the standard topology of $M(n,\C)$. It would be interesting to promote its validity  to the Zariski topology. This would be a closer statement to Arnol'd's original claim.
\end{rem}

\begin{rem}
Theorem \ref{dtd22} implies that the decomposition $M(n,\C)=\coprod_{\bm \la\in\mathscr P(2,n)}\mc M_{\bm \la}$ defines an $\mathscr S$-{\it decomposition} of the space of matrices in the sense of Goresky--MacPherson \cite[Def.\,1.1]{GM88}. Here $\mathscr S$ is a poset structure defined on the set $\mathscr P(2,n)$, with respect a partial order relation that we are going to describe.
\end{rem}

Define the relation $\unlhd$ on the set of bundles in $M(n,\C)$ by
\[\mc M_{\bm \la} \unlhd \mc M_{\bm \nu}\quad \colon\!\!\!\Longleftrightarrow\quad \mc M_{\bm \la}\subseteq \overline{\mc M_{\bm \nu}}.
\]The relation $\unlhd$ defines a partial ordering on bundles, called {\it closure relation}. Such a relation has been extensively studied in \cite{BHM98,DE95,EEK97,EEK99,EJK03,Per04}. A convenient way to visualize\footnote{We invite the reader to use the software StratiGraph, developed at the Ume\r{a} University (Sweden), to visualize the Hasse diagrams in low dimensions. The software is available at the web-page https://www.umu.se/en/research/projects/stratigraph-and-mcs-toolbox/} the closure stratification of bundles is via the Hasse diagram of the relation $\unlhd$: each bundle $\mc M_{\bm \la}$ is represented by a point in the plane (the vertices of the diagram), and one draws an arrow from $v_1$ to $v_2$ if $v_2$ {\it covers} $v_1$ (i.e.\,\,we have $v_1\unlhd v_2$ and there is no $v_3$ such that $v_1\unlhd v_3\unlhd v_2$). See Figure \ref{Hasse}.
\vskip2mm
\begin{figure}\[\xymatrix@R=0.8em@C=0.8em{
16&15&14&13&12&11&10&9&8&7&1\\
&&&\al\al\bt\gm\ar@/_1pc/[dll]&\al^2\bt\bt\ar[l]&&\al\al\bt\bt\ar[ll]&&&&\\
\al\bt\gm\dl&\al^2\bt\gm\ar[l]&\al^3\bt\ar[l]&&\al^2\al\bt\ar[ul]\ar[ll]&&&&\al\al\al\bt\ar[llll]&&\\
&&\al^2\bt^2\ar@/^/[ul]&\al^4\ar[l]\ar[ul]&&\al^3\al\ar[ll]\ar[lu]\ar@/_1pc/[uul]&&\al^2\al^2\ar[ll]\ar@/_1pc/[uul]&&\al^2\al\al\ar[ll]\ar@/_1pc/[ul]&\al\al\al\al\ar[l]
}
\]
\caption{Hasse diagram of the closure stratification of bundles in $M(4,\C)$. The numbers on the top row denote the dimensions of the bundles.}\label{Hasse}
\end{figure}
Introduce the following two types of elementary transformations of $\bm \la$:
\vskip2mm
\noindent{\bf Type I.} $\bm \la\to\bm \la'$, where $\bm \la'$ is obtained by ``merging'' two distinct partitions inside $\bm \la$, say $\msl\la_{i1},\dots,\la_{il_i}\msr$ and $\msl\la_{j1},\dots,\la_{jl_j}\msr$ with $i\neq j$, in the single one 
\[\la_{i1}+\la_{j1}\geq \la_{i2}+\la_{j2} \geq \la_{i3}+\la_{j3}\geq \dots
\]The resulting rough length $|\!|\bm\la'|\!|$ 
equals $|\!|\bm\la|\!|-1$.
\vskip2mm
\noindent{\bf Type II.} $\bm \la\to\bm \la'$, where $\bm \la'$ is obtained as follows. Fix $i\in\{1,\dots,k\}$, and consider the Ferrers diagram corresponding\footnote{In the first column we have $\la_{i1}$ boxes, in the second column we have $\la_{i2}$ boxes, \dots and so on.} to the partition $\la_i=\msl\la_{i1},\dots,\la_{il_i}\msr$ inside $\bm\la$.  
Then move one box rightward one column, or downward one raw, so long as the corresponding partitions $\la_i'$ remains monotonic. The double partition $\bm \la'$ is obtained by $\bm\la$ by replacing $\la_i\mapsto \la_i'$. The resulting rough length $|\!|\bm\la'|\!|$ 
equals $|\!|\bm\la|\!|$.

\begin{thm}[{\cite[Th.\,2.6]{EEK99}}]
Let 
\[\bm \la=\msl \la_{11},\dots, \la_{1l_1};\dots;\la_{k1},\dots,\la_{kl_k}\msr,\quad \bm \la'=\msl\la'_{11},\dots, \la'_{1l_1'};\dots;\la'_{k'1},\dots,\la'_{k'l_k'}\msr,
\]be two double partitions of $n$, with $k'\leq k$. We have $\mc M_{\bm \la'}\unlhd \mc M_{\bm \la}$ if and only if $\bm \la'$ can be obtained from $\bm \la$ via a finite sequence of transformations of type I and II. 

If $\bm\la$ and $\bm \la'$ are related by a single transformation of type I and II, then there is an arrow from $\bm\la'$ to $\bm\la$ in the Hasse diagram of the closure relation $\unlhd$.\qed
\end{thm}
\begin{rem}
Transformations of Type I correspond to coalescences/splittings of eigenvalues.
\end{rem}

The following result gives some insights into the difficulty of the decision procedure for testing the closure relation.
\begin{thm}[{\cite[Th.\,2.7]{EEK99}}]
Deciding whether a bundle is in the closure of another bundle is an NP-complete problem. \qed
\end{thm}

\subsection{Conjugate bundles, and sets $\Mreg,$ $\Mdiag$}\label{secmreg}The set $\mathscr P(n)$ of ordinary partitions of $n$ is equipped with a natural involution $(-)^\vee\colon \mathscr P(n)\to \mathscr P(n)$, called {\it conjugation}.
\vskip2mm
The easiest way to define it is in terms of Ferrers diagrams. Given $\la\in\mathscr P(n)$, with associated Ferrers diagram $F_\la$, the conjugate partition $\la^\vee$ is the one associated with the transposed Ferrers diagram $F_\la^T$ (i.e. the diagram obtained by flipping $F_\la$ along its main diagonal, by turning rows to columns, and vice-versa).
\vskip2mm
This involution $\la\mapsto\la^\vee$ naturally extends to double partitions. Given a double partition $\bm \la=\msl \la_1,\dots, \la_r\msr$ in $\mathscr P(2,n)$, its conjugate $\bm\la^\vee$ is obtained by applying $(-)^\vee$ elementwise, that is
\[\bm \la^\vee=\msl \la_1^\vee,\dots, \la_r^\vee\msr.
\]

We say that two bundles $\mc M_{\bm \la_1}$ and $\mc M_{\bm \la_2}$ are conjugate, if $\bm \la_1$ and $\bm \la_2$ are conjugate double partitions via $(-)^\vee$.
\vskip2mm
Introduce the following two subsets $\mc R,\mc D\subseteq \mathscr P(2,n)$:
let $\mc R$ to be the set of double partitions of $n$ whose fine lengths equal 1, and let $\mc D$ to be the set of double partitions whose fine parts equal $1$.

\begin{lem}
Elements of $\mc R$ are conjugate of elements of $\mc D$, and vice-versa.\qed
\end{lem}

Define the following unions of bundles
\[\Mreg:=\coprod_{\bm \la\in\mc R}\mc M_{\bm \la},\qquad \Mdiag:=\coprod_{\bm \la\in\mc D}\mc M_{\bm \la}.
\]The set $\Mreg$ equals the set of {\it regular matrices}, that is the set of matrices $A\in M(n,\C)$ satisfying one (and hence all) of the following equivalent conditions:
\begin{enumerate}
\item the characteristic polynomial of $A$ equals its minimal polynomial,
\item the centralizer of $A$ in $M(n,\C)$ is of minimal dimension (i.e. it equals $n$),
\item the centralizer of $A$ in $M(n,\C)$ is $\C[A]$.
\end{enumerate}
The set $\Mdiag$ is the set of diagonalizable matrices, with possibly non-simple spectrum.

\section{On the similarity and the Jordan forms of holomorphic matrices}\label{secJ}

\subsection{Global and local holomorphic similarity}
Let $X$ be a complex manifold. Two holomorphic maps $A_1,A_2\colon X\to M(n,\C)$ are said to be 
\begin{itemize}
\item  ({\it globally}) {\it holomorphically similar} on $X$ if there exists a holomorphic map $H\colon X\to GL(n,\C)$ such that $A_1=H^{-1}A_2 H$.
\item {\it locally holomorphically similar at $x_o\in X$} if there exists a neighborhood $U\subseteq X$ of $x_o$ such that $A_1|_U$ and $A_2|_U$ are holomorphically similar on $U$. 
\end{itemize}
Analogue definitions of {\it continuous} or $\mc C^k$-{\it smooth similarity}, with $0\leq k\leq \infty$, can be given, according to the regularity of the matrix-valued function $H$ above.
\vskip2mm

The problem of upgrading local to global holomorphic similarity has been addressed in \cite{Gur88,Lei20,For17}. The following positive results have been obtained, provided that $X$ is a Stein space.

\begin{thm}[\cite{Gur88,Lei20,For17}]\label{thste1}
Let $X$ be a one dimensional Stein manifold, and $A_1,A_2\colon X\to M(n,\C)$ be two holomorphic maps. If $A_1,A_2$ are locally holomorphically similar at each point of $X$, then they are globally holomorphically similar on $X$.\qed
\end{thm}

In the original proof of R.M.\,Guralnick, $X$ is assumed to be a noncompact connected Riemann surface only. Recently, J.\,Leiterer extended Guralnick's result to all one dimensional Stein spaces (not necessarily smooth): his proof is based on the Oka principle for Oka pairs by O.\,Forster and K.J.\,Ramspott \cite{FR66}. An alternative proof of this result also appears in the book \cite{For17}, where it is invoked an alternative Oka principle established in \cite{For03}.
\vskip2mm
For Stein spaces of arbitrary dimensions, we have the following result, which requires stronger assumptions.
\begin{thm}[{\cite[Th.\,1.4]{Lei20}}]\label{thste2}Let $X$ be a Stein space, and $A_1,A_2\colon X\to M(n,\C)$ be two holomorphic maps such that:
\begin{enumerate}
\item $A_1$ and $A_2$ are globally continuously similar on $X$, i.e.\,there exists a continuous map $C\colon X\to M(n,\C)$ such that $A_1=C^{-1}A_2C$,
\item $A_1$ and $A_2$ are locally holomorphically similar at each point of $X$, i.e.\,for each $x_o\in X$ there is a neighborhood $U_o$ of $x_o$, and a holomorphic map $H_o\colon U_o\to GL(n,\C)$ with $A_1=H_o^{-1}A_2 H_o$ on $U_o$,
\item we have $H_o(x_o)=C(x_o)$ for each $x_o\in X$.
\end{enumerate}
Then $A_1$ and $A_2$ are globally holomorphically similar on $X$.\qed
\end{thm}

Conditions (1) and (2) alone do not imply global holomorphic similarity. For a counterexample, see \cite[Th.\,8.2]{Lei20}.

\subsection{Three criteria for local holomorphic similarity}Let $X$ be a complex manifold, $A_1,A_2\colon X\to M(n,\C)$ be two holomorphic maps, $x_o\in X$, and $\Phi\in M(n,\C)$ such that $\Phi A_1(x_o)=A_2(x_o)\Phi$.
\vskip2mm
Below are some criteria on $X, A_1,A_2,\Phi$ which will allow to infer local holomorphic similarity of $A_1$ and $A_2$ on a neighborhood of $x_o$. Following Leiterer, we will call them {\it Wasow's}, {\it Smith's}, and {\it Spallek's criterion} respectively.
\vskip2mm
\noindent{\bf Wasow's criterion: }
The dimension of the complex vector space $$\{\Theta\in M(n,\C)\colon \Theta A_1(x)=A_2(x)\Theta\}$$ is constant for $x$ in some neighborhood of $x_o$.
\vskip2mm
\noindent{\bf Smith's criterion: }
The space $X$ is one dimensional, and there exist a neighborhood $V_o$ of $x_o$ and a continuous map $C_o\colon V_o\to M(n,\C)$ such that $C_oA_1=A_2C_o$ on $V_o$, and $C_o(x_o)=\Phi$. 
\vskip2mm
\noindent{\bf Spallek's criterion: }
There exist a neighborhood $V_o$ of $x_o$ and a smooth map $T_o\colon V_o\to M(n,\C)$ such that $T_oA_1=A_2T_o$ on $V_o$, and $T_o(x_o)=\Phi$.

\begin{thm}[{\cite[Th.\,1.5]{Lei20}}]
\label{thcrit}If one of the criteria above holds, then there exists a neighborhood $U_o$ of $x_o$ and a holomorphic map $H_o\colon U_o\to M(n,\C)$ such that $H_oA_1=A_2H_o$ on $U_o$, and $H_o(x_o)=\Phi$. In particular, if $\Phi$ is invertible, then $A_1$ and $A_2$ are locally holomorphically similar. \qed
\end{thm}
\begin{rem}
The names of the criteria are justified as follows. In \cite{Was62}, W.\,Wasow formulated the first criterion and proved the statement of Theorem \ref{thcrit} under the unnecessary assumption that $X$ is a domain in $\C$. The proof of Theorem \ref{thcrit} under the Smith's criterion is based on applications of the Smith factorization theorem, see \cite[Ch.\,III, Sec.\,8]{Jac75} \cite[Th.\,4.3.1]{GL09}. On the other hand, the proof of Theorem \ref{thcrit} under the Spallek's criterion follows from a special case of a result of K.\,Spallek, see \cite[Satz 5.4]{Spa65} \cite[Introduction]{Spa67}.
\end{rem}

Spallek's criterion and Theorem \ref{thste2} imply the following result.
\begin{cor}
Let $X$ be a Stein manifold. Let $A_1,A_2\colon X\to M(n,\C)$ be two holomorphic maps, globally $\mc C^\infty$-smoothly similar on $X$. Then $A_1$ and $A_2$ are globally holomorphically similar on $X$. \qed
\end{cor}

\begin{rem}
The statement of this corollary is optimal: the $\mc C^\infty$-smoothness condition cannot be replaced by a $\mc C^k$-smoothness with $k<\infty$. See \cite[Th.\,8.2]{Lei20}.
\end{rem}

Smith's criterion allows to strengthen Theorem \ref{thste1} as follows.
\begin{cor}\label{corleit}
Let $X$ be a one dimensional Stein manifold, and let $A_1,A_2\colon X\to M(n,\C)$ be two holomorphic maps. If $A_1$ and $A_2$ are locally continuously similar at each point $x\in X$, then they are globally holomorphically similar on $X$. \qed
\end{cor}

\subsection{Holomorphically Jordanizable matrices }\label{secholjor}Let $X$ be a complex manifold, 
and $A\colon X\to M(n,\C)$ be a holomorphic map.
\vskip2mm
We say that $A$ is {\it locally holomoprhically Jordanizable} (or {\it similar to a Jordan matrix}) at $x_o\in X$ if there exists a neighborhood $U\subseteq X$ of $x_o$, and a holomorphic map $S\colon U\to GL(n,\C)$ such that
\beq
\label{lhj1}
J(x)=S(x)^{-1}A(x)S(x),\quad x\in U,
\eeq
is a Jordan matrix for each $x\in U$. If one can take $U=X$, then we say that $A$ is {\it globally} holomorphically Jordanizable.
\vskip2mm
Given $\Phi\in M(n,\C)$, denote by $\si(\Phi)$ the spectrum of $\Phi$, i.e.\,the set of its eigenvalues. Furthermore, for $\la\in\si(\Phi)$ and $k=1,\dots, n$, denote by $\vartheta_k(\Phi,\la)$ the number of Jordan blocks of size $k$ and eigenvalue $\la$ in the Jordan form of $\Phi$, and set
\[\vartheta_k(\Phi):=\sum_{\la\in\si(\Phi)}\vartheta_k(\Phi,\la).
\]
\begin{lem}\label{lemobv}If $A\colon X\to M(n,\C)$ is locally holomorphically Jordanizable and \eqref{lhj1} holds, then we have the following necessary conditions:
\begin{enumerate}
\item The Jordan matrix $J(x)$ is holomorphic on some open set $U\subseteq X$,
\item there exist some holomorphic functions $\la_1,\dots,\la_m\colon U\to \C$ such that $\si(A(x))=\{\la_1(x),\dots, \la_m(x)\}$ for any $x\in U$,
\item for each $k=1,\dots,n$, the functions $U\ni x\mapsto \vartheta_k(A(x))$ are constant.
\end{enumerate}
\end{lem}
\proof
Condition (1) directly follows from \eqref{lhj1}. Condition (1) implies both 
(2) and (3). 
\endproof
Notice that 
condition (2) implies that ${\rm card\,}\si(A(x))\leq m$ for $x\in U$. In particular, the sign < holds in case of coalescences of some of the eigenvalues $\la_1,\dots, \la_m$.
\vskip2mm
Neither condition (1), nor conditions (2) and (3) together, are sufficient for the holomorphic similarity of $A$.

\begin{cex}\label{cex1}
Consider the unit disc $\mathbb D=\{z\in\C\colon |z|<1\}$ in $\C$, and let $A$ and $J$ be the holomorphic matrices given by
\[A(z)=\begin{pmatrix}
z&1&0\\
0&z^2&z\\
0&0&z^2
\end{pmatrix},\quad 
J(z)=\begin{pmatrix}
z&0&0\\
0&z^2&1\\
0&0&z^2
\end{pmatrix},\quad z\in\mathbb D.
\]
It is easy to see that $A(z)$ is similar to $J(z)$ for any $z\in\mathbb D$. However, $A(z)$ is not holomorphically similar to $J(z)$. Assume there exists a holomorphic matrix $S(z)=(s_{ij}(z))_{i,j=1}^3$ such that $S(z)J(z)=A(z)S(z)$ for $z\in\mathbb D$. We have
\begin{align*}
zs_{11}+s_{21}=zs_{11},&&z^2s_{21}+zs_{31}=zs_{21},\\
z^2s_{22}+zs_{32}=z^2s_{22},&&z^2s_{23}+zs_{33}=s_{22}+z^2s_{23},
\end{align*}
for $z\in\mathbb D$. Hence, $s_{21}=s_{31}=s_{32}=0$, and $s_{22}(0)=0\cdot s_{33}(0)=0$, so that $S(0)$ is not invertible.
\end{cex}

\begin{cex}\label{cex2}
Conditions (2) and (3) together do not imply condition (1) above. Consider the disc $\mathbb D=\{z\in\C\colon |z|<\frac{1}{2}\}$ in $\C$, and the holomorphic map $A\colon \mathbb D\to M(4,\C)$ defined by
\[A(z)=\begin{pmatrix}
z&1&0&0\\
0&-z&0&0\\
0&0&1+z&z\\
0&0&0&1+z
\end{pmatrix},\quad z\in\mathbb D.
\]The eigenvalues of $A(z)$ are given by the holomorphic functions $\la_1(z)=z,\,\la_2(z)=-z,\,\la_3(z)=1+z$. Furthermore, we have $\vartheta_1(A(z))=2$ and $\vartheta_2(A(z))=1$ for each $z\in\mathbb D$. The matrix $A$, however, does not admit a holomorphic Jordan form.
\end{cex}

\subsection{Coalescing points}Let $X$ be a complex manifold, 
and $A\colon X\to M(n,\C)$ be a holomorphic map. 
\begin{defn}
A point $x_o\in X$ is called a {\it coalescing point} of the eigenvalues of $A$ if for any neighborhood $U$ of $x_o$ there exists $x\in U$ such that ${\rm card\,}\si(A(x))> {\rm card\,}\si(A(x_o)).$ We denote by coal$(A)$ the set of coalescing points of eigenvalues of $A$.
\end{defn}

For a proof of the following well-known result, see e.g.\,\,\cite{Lei17}.
\begin{prop}\label{propcoal}
Let $A\colon X\to M(n,\C)$ be a holomorphic map, $x_o\in X$, $\si(A(x_o))=\{\la_{o,1},\dots, \la_{o,m}\}$, and let $n_j$ be the algebraic multiplicity (i.e.\,the order as a zero of the characteristic polynomial of $A(x_o)$) of $\la_{o,j}$ for $j=1,\dots,m$.

The following conditions are equivalent:
\begin{enumerate}
\item $x_o\notin {\rm coal}(A(x_o))$;
\item there exists a neighborhood $U$ of $x_o$, and uniquely determined holomorphic functions $\la_1,\dots,\la_m\colon U\to \C$ such that
\begin{itemize}
\item $\la_j(x_o)=\la_{o,j}$ for $j=1,\dots,m$,
\item $\si(A(x))=\{\la_1(x),\dots, \la_m(x)\}$ for each $x\in U$,
\item the eigenvalue $\la_j(x)$ of $A(x)$ has algebraic multiplicity $n_j$ for each $x\in U$.\qed
\end{itemize}
\end{enumerate}
\end{prop}

\begin{thm}
If the set ${\rm coal}(A)$ is non-empty, then it is a 
closed analytic subset of $X$ of codimension 1. \qed
\end{thm}
\begin{rem}
There exist many sources in the literature for a proof of this result, e.g.\,\,see \cite{Bau74,Bau85}\cite[Ch.\,III, Th.\,4.3 and 4.6]{FG02}. In the recent preprint \cite[Th.\,4.3]{Lei17}, J.\,Leiterer considered also general complex spaces $X$ (i.e.\,\,by allowing singularities). Leiterer's result provides finer estimates: if $X$ is irreducible, and ${\rm coal}(A)\neq \emptyset$, then there exist finitely many holomoprhic functions $h_1,\dots,h_\ell\colon X\to \C$ such that
\[{\rm coal}(A)=\{x\in X\colon h_j(x)=0,\,j=1,\dots, \ell\}.
\]Each $h_j$ is a finite sum of products of the entries of $A(x)$. Moreover, we have 
\[|h_j(x)|\leq (2n)^{6n^2}|\!|A(x)|\!|^{2n^2},\quad x\in X,\quad1\leq j\leq \ell.
\]
\end{rem}

\subsection{The gap topology}
In what follows, we denote by $|\!| f|\!|$ the operator norm of any linear map $f\colon \C^{n}\to \C^m$, which is defined by $|\!|f|\!|:=\sup_{|\!|x|\!|=1}|\!|f(x)|\!|$. Here the spaces $\C^n,\C^m$ are intended to be equipped with the standard hermitian metric.
\vskip2mm
Denote by $\mc G_n$ the set of all $\C$-vector subspaces of $\C^n$. 
Given $L_1,L_2\in\mc G_n$, denote by $\Pi_1,\Pi_2\colon \C^n\to\C^n$ the orthogonal projections onto $L_1$ and $L_2$ respectively. We define the {\it gap distance} between $L_1$ and $L_2$ as
\[\Theta(L_1,L_2):=|\!|\Pi_1-\Pi_2|\!|.
\]

\begin{prop}[{\cite[Ch.\,XIII]{GLR06}}]
The gap distance $\Theta$ defines a metric on $\mc G_n$. 
Moreover, we have 
\begin{enumerate}
\item $\Theta(L_1,L_2)\leq 1$ for each $L_1,L_2\in\mc G_n$,
\item $\Theta(L_1,L_2)<1$ only if $\dim L_1=\dim L_2$. \qed
\end{enumerate}
\end{prop}

Denote by $G(k,n)$ the Grassmannian of complex $k$-dimensional subspaces in $\C^n$. The space $G(k,n)$ can be defined as the topological quotient $U(n)/U(k)\times U(n-k)$, with respect to the inclusion $U(k)\times U(n-k)\hookrightarrow U(n)$. In particular, it follows that $G(k,n)$ is compact and connected.
\begin{cor}\label{gngrass}
The connected components of the metric space $(\mc G_n,\Theta)$ are homeomorphic to the Grassmannians of complex subspaces in $\C^n$, i.e.
$\mc G_n\cong \coprod_{k=0}^nG(k,n).
$ In particular, the metric space $(\mc G_n,\Theta)$ is complete.
\end{cor}
\proof By the previous proposition, the function $L\mapsto \dim L$ is locally constant on $\mc G_n$, and consequently constant on the connected components of $\mc G_n$. 
Denote $\mc G_{n,k}$ the connected component of $\mc G_n$ whose points are $k$-dimensional subspaces of $\C^n$.
We claim that the identity map $G(k,n)\to \mc G_{n,k}$ is continuous, for each $k=0,\dots,n$. To see this, it is sufficient to show the continuity of the map $f\colon U(n)\to \mc G_{n,k}$ defined by 
\[A\mapsto \text{$\C$-span of the first $k$ columns of $A$.}
\]Let $A_1,A_2\in U(n)$, and denote by $\widetilde A_1,\widetilde A_2$ the $n\times k$ matrices obtained by the first $k$ columns of $A_1$ and $A_2$, respectively. Let $\Pi_1,\Pi_2\in \End(\C^n)$ be the orthogonal projections onto $f(A_1)$ and $f(A_2)$, respectively. With respect to the standard basis of $\C^n$, the matrices representing $\Pi_1$ and $\Pi_2$ equal $\widetilde A_1 \widetilde A_1^*$ and $\widetilde A_2 \widetilde A_2^*$. 
We have
\begin{align*}
\Theta(f(A_1),f(A_2))&=|\!|\Pi_1-\Pi_2|\!|=|\!|\widetilde A_1 \widetilde A_1^*-\widetilde A_2 \widetilde A_2^*|\!|=|\!|\widetilde A_1(\widetilde A_1^*-\widetilde A_2^*)+(\widetilde A_1-\widetilde A_2)\widetilde A_2^*|\!|\\
&\leq |\!|\widetilde A_1|\!| |\!|\widetilde A_1^*-\widetilde A_2^*|\!|+ |\!|\widetilde A_1-\widetilde A_2|\!| |\!|\widetilde A_2^*|\!|\\
&\leq |\!|A_1|\!| |\!|A_1^*-A_2^*|\!|+ |\!|A_1-A_2|\!| |\!|A_2^*|\!| = 2 |\!|A_1-A_2|\!|. 
\end{align*}
This proves that the identity map $G(k,n)\to \mc G_{n,k}$ is continuous. Moreover it is also closed, since $G(k,n)$ is compact and $\mc G_n$ is Hausdorff.
\endproof

\subsection{Holomorphic families of subspaces} The topological space $\mc G_n$ can be made into a complex manifold in a natural way, due to Corollary \ref{gngrass}. 
\vskip2mm
We call a {\it holomorphic family of subspaces of $\C^n$}, parametrized by a complex manifold $X$, any holomorphic map $f\colon X\to\mc G_n$. If $X$ is connected, then $f$ takes values in a complex Grassmannian $G(k,n)$ for some $k=0,\dots,n$.
\begin{rem}\label{equivdef}
There is a 1-1 correspondence between the following data:
\begin{enumerate}
\item a holomorphic family of subspaces $f\colon X\to G(k,n)$;
\item a rank $k$ holomorphic subbundle of the trivial bundle $\underline{\C^n}:=X\times \C^n$;
\item a continuous map $f\colon X\to \mc G_n$ such that, for any $x_o\in X$ there exist a neighborhood $U$ of $x_o$ and holomorphic maps $v_1,\dots, v_k\colon U\to\C^n$ such that $v_1(z),\dots, v_k(z)$ are linearly independent and $f(z)={\rm span}\langle v_1(z),\dots, v_k(z)\rangle$ for each $z\in U$.
\end{enumerate}
\end{rem}

\begin{prop}\label{propcont}
Let $X$ be a locally compact metric space, and $T\colon X\to M(n,\C)$ a continuous map. The following conditions are equivalent:
\begin{enumerate}
\item $\dim\ker T(x)$ is locally constant on $X$;
\item $\dim\ims\,T(x)$ is locally constant on $X$;
\item the map $f\colon X\to \mc G_n$, $x\mapsto \ker T(x)$, is continuous;
\item the map $f\colon X\to \mc G_n$, $x\mapsto \ims\,T(x)$, is continuous.
\end{enumerate}
If moreover $X$ is a complex manifold and $T$ is holomorphic, then the conditions above are equivalent to the following ones:
\begin{enumerate}
\item[(5)] the map $f\colon X\to \mc G_n$, $x\mapsto \ker T(x)$, is holomorphic;
\item[(6)] the map $f\colon X\to \mc G_n$, $x\mapsto \ims\,T(x)$, is holomorphic.
\end{enumerate}
\end{prop}
\proof
The only non-trivial statements are $(1)\Rightarrow (3)\Rightarrow (5)$ and $(2)\Rightarrow (4)\Rightarrow (6)$. The implications $(1)\Rightarrow (3)$ and $(2)\Rightarrow (4)$ can be easily proved by adapting the argument of \cite[Prop.\,13.6.1]{GLR06}. The implications $(3)\Rightarrow (5)$ and $(4)\Rightarrow (6)$ follow from a result of Ph.G.A.\,Thijsse, \cite[Th.\,3.1]{Thi78}, proved in the more general case of families of complemented subspaces of Banach spaces. See \cite[Prop.\,5 and Appendix]{Jan88} for the complete argument, and the reproduction of the proof of Thijsse's result. See also the argument of \cite[Prop.\,18.1.2]{GLR06} for the case $\dim X=1$. 
\endproof

Given a holomorphic map $T\colon X\to M(n,\C)$, define $n_0:=\min_{x\in X}\{\dim\ker T(x)\}$. The set $\Ups(T)$ of \emph{jump points} of $T$ is the set
\[\Ups(T):=\{x\in X\colon \dim\ker T(x)>n_0\}.
\]
\begin{prop}[{\cite[Satz 1.1]{Kab76}}]
The set $\Ups(T)$ is an analytic subset of $X$.\qed
\end{prop}
For $x\in \Ups(T)$, the space $\ker T(x)$ is ``too big'' if compared to $\ker T(y)$ for $y\notin\Ups(T)$ near $x$. We can however introduce a suitable replacement for the space $\ker T(x)$. 
\vskip2mm
With any holomorphic map $T\colon X\to M(n,\C)$, there is a naturally associated morphism of $\mathscr O_X$-modules (for simplicity denoted by the same symbol)
\[T\colon\mathscr O_X^{\oplus n}\to \mathscr O_X^{\oplus n},
\]where $\mathscr O_X$ denote the structure sheaf of $X$. For any $x\in X$ denote by {\calligra\large ker}$\,\,T_x$ the stalk of the kernel sheaf {\calligra\large ker}$\,\,T$. Define $\mc K[T;x]$ to be the space of values at $x$ of germs of $\C^n$-valued holomorphic functions in {\calligra\large ker}$\,\,T_x$, i.e.
\[
\mc K[T;x]=\{v\in \ker T(x)\,|\,\exists f_x\in\text{{\calligra\large ker}$\,\,T_x$},\,\,f_x(x)=v\}
\]
\begin{thm}[\cite{Kab76,Kab12}]\label{kabth}$\quad$
\begin{enumerate}
\item We have 
$\mc K[T;x]\subseteq \ker T(x)$ for any $x\in X$.
\item We have $
\mc K[T;x]= \ker T(x)$ for $x\notin \Ups(T)$.
\item There exists an analytic subset $\Si\subseteq X$ of codimension at least 2 such that $f\colon X\setminus \Si\to \mc G_n$, $x\mapsto \mc K[T;x]$, 
is a holomorphic family of subspaces.\qed
\end{enumerate}
\end{thm}

\subsection{A generalization of a theorem of Thijsse and Wasow} As shown by Counterexamples \ref{cex1} and \ref{cex2}, conditions (1), (2), (3) of Lemma \ref{lemobv} are not sufficient to infer the existence of a holomorphic Jordan form. In this section, we find a further condition which, jointly with condition (1), will ensure the locally holomoprhically Jordanizability of a matrix. 
\vskip2mm
\begin{thm}\label{th1}
Let $X$ be a complex manifold, and $A\colon X\to M(n,\C)$ be a holomorphic map. The matrix $A$ is locally holomorphically Jordanizable at a point $x_o\in X$ if and only if the following conditions are satisfied:
\begin{enumerate}
\item there exists a neighborhood $U$ of $x_o$ and a holomorphic map $J\colon U\to M(n,\C)$ such that $J(x)$ is a Jordan form of $A(x)$ for each $x\in U$; in particular, there exist local holomorphic functions $\la_1,\dots,\la_r\colon U\to \C$ such that $\si(A(x))=\{\la_1(x),\dots, \la_r(x)\}$ for each $x\in U$;
\item for each $i=1,\dots, r$ the limits of the generalized eigenspaces 
\[\lim_{\substack{z\to x_o\\z\notin {\rm coal}(A)}}\ker \left(A(z)-\la_i(z){\rm Id}\right)^n
\]exist in the gap topology of $\mc G_n$;
\item we have
\[\C^n=\bigoplus_{i=1}^r\lim_{\substack{z\to x_o\\z\notin {\rm coal}(A)}}\ker \left(A(z)-\la_i(z){\rm Id}\right)^n.
\]
\end{enumerate}
Moreover, conditions $(2)$ and $(3)$ are respectively equivalent to the following ones:
\begin{enumerate}
\item[(2bis)] the function $x\mapsto \dim \mc K[(A-\la_i{\rm Id})^n;x]$ 
is continuous at $x_o$, for any $i=1,\dots, r$;
\item[(3bis)] we have $
\C^n=\bigoplus_{i=1}^r \mc K[(A-\la_i{\rm Id})^n;{x_o}]$. 
\end{enumerate}
\end{thm}

\begin{rem} If the limit of point (2) exists it necessarily equals $\mc K[(A-\la_i{\rm Id})^n;{x_o}]$. 
This follows form Proposition \ref{propcont} and Theorem \ref{kabth}. This implies the equivalences of (2) and (2bis), and of (3) and (3bis).
Moreover, notice that conditions (2), (2bis), (3), and (3bis) are trivially satisfied if $x_o\notin {\rm coal}(A)$. This follows from Propositions \ref{propcoal} and \ref{propcont}.
\end{rem}

Before proving the theorem, we consider a simpler case, namely that of a holomorphic Jordan form with one eigenvalue only. 

\begin{lem}\label{lemj}
Let $T\colon X\to M(n,\C)$ be a holomorphic map. Assume that
\begin{enumerate}
\item there exists a holomorphic map $J\colon X\to M(n,\C)$ such that $J(x)$ is a Jordan form of $T(x)$ for each $x\in X$,
\item $T$ has a unique holomorphic eigenvalue function $\la\colon X\to \C$.
\end{enumerate}
Then $T$ is holomorphically Jordanizable on any domain of $X$ biholomorphic to a polydisc. 
\end{lem}

\proof For any $z\in X$, set $M(z):=T(z)-\la(z){\rm Id}$, and $N_j(z):=\ker M(z)^j$ for any $j=1,2,3,\dots$. We have the tower of subspaces
\[N_1(z)\subseteq N_2(z)\subseteq N_3(z)\subseteq\dots\subseteq N_{n-1}(z)\subseteq N_{n}(z)=N_{n+1}(z)=\dots.
\]

First, we show that we have well-defined holomorphic maps $N_j\colon X\to\mc G_n$ for any $j$. 
\vskip2mm
For each $j$, let us introduce the analytic subsets $\Ups_j\subseteq X$ of jump points of $M^j$. We claim that $\Ups_j=\emptyset$ for any $j$. 

For $j\geq n$ the statement is obvious. We have the following facts:\newline
\noindent (i) since the matrix $T$ has a single eigenvalue $\la$, we have
\[\vartheta_k(T(z))=\vartheta_k(T(z),\la(z)),\quad z\in X,\quad k=1,\dots,n;
\]
\noindent (ii) since $J(z)$ is holomorphic, the function $z\mapsto \vartheta_k(T(z))$ is constant on $X$ for any $k$; \newline
\noindent (iii) we have 
\[\vartheta_k(T(z),\la(z))=\dim \frac{N_{k}(z)}{N_{k-1}(z)}-\dim \frac{N_{k+1}(z)}{N_k(z)},\quad k=1,\dots, n,\quad N_0:=0.
\]From (i), (ii), (iii) for $k=n$, and the fact that $\Ups_n=\emptyset$, we deduce that $\Ups_{n-1}=\emptyset$. Then, by applying (i), (ii), (iii) for $k=n-1$, we deduce that $\Ups_{n-2}=\emptyset$. By iteration of this argument, one proves that all the sets $\Ups_j$ are empty. So $N_j\colon X\to\mc G_n$ are holomorphic for any $j$. In particular, by Remark \ref{equivdef}, each $N_j$ can be seen as a holomorphic subbundle of $\underline{\C^n}:=X\times \C^n$.
\vskip2mm
Let $\Dl\subseteq X$ a domain biholomorphic to a polydisc. In particular, $\Dl$ is a Stein manifold. For any $j\geq 1$, there exists a holomorphic vector subbundle $V_j\to \Dl$ of $\underline{\C^n}|_{\Dl}$ such that $N_j|_\Dl\oplus V_j= N_{j+1}|_{\Dl}$. This follows from a general result of Shubin \cite[Th.\,1]{Shu70}. 
Since $\Dl$ is Stein, the topological and analytical classifications of vector bundles coincide, by the Oka-Grauert principle. Consequently, all vector bundles $N_j|_{\Dl}$ and $V_j$, with $j\geq 1$, are trivial. 
\vskip2mm
Let $n_j$ be the dimension of the subspaces $N_j(z)$ for any $j\geq 1$, and set $\ell:=\min\{j\colon n_j=n_{j+1}\}$. Since $V_{\ell-1}$ is trivial, there exist a global basis of sections, i.e. holomorphic maps $v_1,\dots, v_{n_\ell-n_{\ell-1}}\colon \Dl\to \C^n$ such that
\[N_{\ell}(z)=N_{\ell-1}(z)\oplus{\rm span}\langle v_1(z),\dots, v_{n_\ell-n_{\ell-1}}(z)\rangle,\quad z\in X.
\]The tuple 
\[v_1(z),\dots,\,\,v_{n_\ell-n_{\ell-1}}(z),\,\,M(z)v_1(z),\dots,\,\,M(z)v_{n_\ell-n_{\ell-1}}(z)
\]is easily seen to be linearly independent. Hence, by the triviality of $V_{\ell-2}$, there exist vector-valued holomorphic functions $$v_{n_\ell-n_{\ell-1}+1},\dots, v_{n_{\ell-1}-n_{\ell-2}}\colon X\to \C^n,$$ such that
\begin{multline*}
N_{\ell-1}(z)=N_{\ell-2}\oplus\, {\rm span}\langle M(z)v_1(z),\dots,\,\,M(z)v_{n_\ell-n_{\ell-1}}(z)\rangle\\
\oplus\, {\rm span}\langle v_{n_\ell-n_{\ell-1}+1}(z),\dots, v_{n_{\ell-1}-n_{\ell-2}}(z)\rangle,\quad z\in X.
\end{multline*}
Proceeding in this way, by applying the standard construction, a family of holomorphic Jordan bases of $N_\ell$ is obtained.
\endproof

\proof[Proof of Theorem \ref{th1}]
The necessity of condition (1) is clear. The limits 
$$\lim_{\substack{z\to x_o\\z\notin {\rm coal}(A)}}\ker \left(J(z)-\la_i(z){\rm Id}\right)^n,\quad i=1,\dots,r,
$$satisfy conditions (2) and (3). Moreover, if $S\colon U\to GL(n,\C)$ is such that $J=S^{-1}AS$, we have 
\[\lim_{\substack{z\to x_o\\z\notin {\rm coal}(A)}}\ker \left(A(z)-\la_i(z){\rm Id}\right)^n=\lim_{\substack{z\to x_o\\z\notin {\rm coal}(A)}}S(z)\ker \left(J(z)-\la_i(z){\rm Id}\right)^n,
\]and conditions (2) and (3) are satisfied as well. 
\vskip2mm
Let us prove the sufficiency of conditions (1), (2), and (3). If condition (1) is satisfied, there exist holomorphic families of subspaces $L_i\colon X\setminus{\rm coal}(A)\to\mc G_n$ defined by $L_i(z):=\ker \left(A(z)-\la_i(z){\rm Id}\right)^n$ for $i=1,\dots, r$. By Proposition \ref{propcont} and Kaballo's Theorem \ref{kabth}, these families can be prolonged to holomorphic families $L_i\colon X\setminus \Si_i\to \mc G_n$, defined by
\[x\mapsto \mc K[(A-\la_i{\rm Id})^n;x], 
\]on the complement of analytic subspaces $\Si_i\subseteq X$ of codimension $\geq 2$. By condition (2), we necessarily have $x_o\notin \bigcup_{i=1}^r\Si_i$. 

Hence, there exist 
\begin{itemize}
\item a sufficiently small neighborhood $U$ of $x_o$ on which the functions \[n_i(x)=\dim \mc K[(A-\la_i{\rm Id})^n;x], 
\quad i=1,\dots,r,\] are constant, 
\item holomorphic functions $W_i\colon U\to GL(n_i,\C)$, for $i=1,\dots,r$, with $W_i(x_o)={\rm Id}$,
\end{itemize}
such that 
\[L_i(z)=W_i(z)L_i^o,\quad L_i^o\equiv \mc K[(A-\la_i{\rm Id})^n;{x_o}], 
\quad i=1,\dots,r.
\]
Moreover, by condition (3), up to shrinking $U$, we may assume that $\C^n= L_1(z)\oplus\dots\oplus L_r(z)$ for all $z\in U$. Define $W\colon U\to GL(n,\C)$ as the direct sum $W=W_1\oplus\dots\oplus W_r$. We have 
\[W^{-1}(z)A(z)W(z)=\bigoplus_{i=1}^r\widetilde A_i(z),\quad z\in U,
\]where each matrix $\widetilde A_i(z)$ has a unique eigenvalue $\la_i(z)$, for $i=1,\dots,r$. Moreover, for any $i=1,\dots,r$, there exist holomorphic matrix-valued functions $\widetilde J_i\colon U\to M(n,\C)$ which are the Jordan forms of $\widetilde A_i(z)$, by condition (1). 
The result follows by applying Lemma \ref{lemj} to each matrix $\widetilde A_i$, with $i=1,\dots, r$.
\vskip2mm
This completes the proof.
\endproof

\begin{cor}[{\cite[Main Theorem]{Thi85} \cite[Th.\,12.2-2]{Was85}}]\label{corTW}
Let $X$ be a one dimensional Stein manifold, and $A\colon X\to M(n,\C)$ be a holomorphic map. The matrix $A$ is \emph{globally} holomorphically Jordanizable if and only if the following conditions are satisfied:
\begin{enumerate}
\item there exists a holomorphic map $J\colon X\to M(n,\C)$ such that $J(x)$ is a Jordan form of $A(x)$ for each $x\in X$; in particular, there exist holomorphic functions $\la_1,\dots,\la_r\colon X\to \C$ such that $\si(A(x))=\{\la_1(x),\dots, \la_r(x)\}$ for each $x\in X$;
\item we have $
\C^n=\bigoplus_{i=1}^r \mc K[(A-\la_i{\rm Id})^n;{x}]$ 
for each $x\in X$.
\end{enumerate}
\end{cor}
\proof
The result follows from Corollary \ref{corleit}, and Theorem \ref{th1}. Notice that conditions (2) and (2bis) of Theorem \ref{th1} 
are trivially satisfied in the one dimensional case, due to Kaballo's Theorem \ref{kabth}.
\endproof
\vskip2mm

\begin{ex}
Consider the holomorphic map $A\colon \C^2\to M(3,\C)$ defined by $A(\bm x)=\begin{pmatrix}
x_1&0&x_2\\
0&x_1&x_2\\
0&0&0
\end{pmatrix}$. There exists holomorphic eigenvalues functions $\la_1(\bm x)=x_1$ and $\la_2(\bm x)=0$. We have coal$(A)=\{x_1=0\}$, and the only points at which $A(\bm x)$ is not diagonalizable are ${\rm coal}(A)\setminus\{(0,0)\}$. 
The only point $\bm x_o$ at which one of the limits of condition (2) of Theorem \ref{th1} does not exist is $\bm x_o=(0,0)$. For $\bm x\notin {\rm coal}(A)$, we have 
\[\ker \left(A(\bm x)-\la_2(\bm x){\rm Id}\right)^3=\mc K[\left(A-\la_2{\rm Id}\right)^3;\bm x]
={\rm span}\big\langle
(-x_2,-x_2,x_1)^T
\big\rangle,
\]whose limit $\bm x\to (0,0)$ does not exist. Notice that at points $\bm x$ in coal$(A)\setminus\{(0,0)\}$ condition (2) is satisfied, but not condition (3). We have indeed 
\begin{multline*}
\mc K[\left(A-\la_1{\rm Id}\right)^3;\bm x]={\rm span}\big\langle
(1,0,0)^T,(0,1,0)^T
\big\rangle\\
\supseteq\quad
\mc K[\left(A-\la_2{\rm Id}\right)^3;\bm x]={\rm span}\big\langle
(-x_2,-x_2,0)^T
\big\rangle.
\end{multline*}\end{ex}
\subsection{Holomorphic Jordanization and bundles of matrices}
Given two multisets $\mc X_1=(X_1,m_1)$ and $\mc X_2=(X_2,m_2)$, we define their {\it multiunion} $\mc X_1\vee \mc X_2$ as the multiset $(X_1\cup X_2,m)$, where the multiplicity function $m\colon X_1\cup X_2\to\N^*$ is defined as follows 
\begin{empheq}[left={m(x)}{:=}\empheqlbrace]{align*}
m_1(x),\quad &\text{if }x\in X_1\setminus X_2,\\
m_2(x),\quad &\text{if }x\in X_2\setminus X_1,\\
m_1(x)+m_2(x),\quad &\text{if }x\in X_1\cap X_2.
\end{empheq}
We have a natural ``forgetful'' surjection $\xi\colon\mathscr P(2,n)\to\mathscr P(n)$, defined by
\[\msl \la_1,\dots,\la_r\msr\mapsto \la_1\vee\dots\vee\la_r.
\]For each $\la\in\mathscr P(n)$, define
\[\mc F_\la:=\coprod_{\bm \nu\in\xi^{-1}(\la)}\mc M_{\bm \nu}.
\]
\begin{thm}\label{thjb1}
Let $X$ be a connected complex manifold, and $A\colon X\to M(n,\C)$ be a holomorphic map, locally holomorphically Jordanizable at each point of $X$. Then we have $A(X)\subseteq \mc F_\la$ for some $\la\in\mathscr P(n)$.
\end{thm}
\proof
For any $k=1,\dots, n$, the map $x\mapsto \vartheta_k(A(x))$ is constant on $X$, by connectedness of $X$ and Lemma \ref{lemobv}. Assume $x_1,x_2\in X$ to be such that $A(x_1)\in\mc M_{\bm \la_1}$ and $A(x_1)\in\mc M_{\bm \la_2}$. We have $\vartheta_k(A(x_1))=\vartheta_k(A(x_2))$ for any $k=1,\dots, n$. This means that the multiplicities of $k$ as an element of both $\xi(\bm \la_1)$, and $\xi(\bm \la_2)$ resp., are all equal for any $k=1,\dots, n$. Hence $\xi(\bm \la_1)=\xi(\bm \la_2)$.
\endproof

Although the inverse statement is not true (as shown, e.g., by Counterexample \ref{cex1}), we have the following result.

\begin{thm}\label{thjb2}
Assume $A\colon X\to M(n,\C)$ is holomorphic with $A(X)\subseteq \mc M_{\bm \la}$ for some $\bm \la\in\mathscr P(2,n)$. Then $A$ is locally holomorphically Jordanizable at each point of $X$.
\end{thm}

\proof
Since $A(X)\subseteq \mc M_{\bm \la}$, we have ${\rm coal}(A)=\emptyset$. Hence conditions (2) and (3) of Theorem \ref{th1} are trivially satisfied. Also condition (1) holds true, by Proposition \ref{propcoal}.
\endproof

\section{Universality of integrable deformations of solutions of RHB problems}\label{secRHB}
\subsection{Riemann--Hilbert--Birkhoff problems} Consider a disc $D$ in $\Pb^1$, centered at $z=\infty$. Given a holomorphic vector bundle $E^o$ on $D$, equipped with a meromorphic connection $\nabla^o$ admitting a pole at $z=\infty$, the Riemann--Hilbert--Birkhoff (RHB) problem is the following:

\begin{prob}\label{RHB1}
Does there exist a {\it trivial} vector bundle $E^o$ on $\Pb^1$ equipped with a meromorphic connection $\nabla^o$, restricting to the given data on $D$, and with a further {\it logarithmic} pole only at $z=0$? 
\end{prob}

Assume that the pole at $z=\infty$ is of order 2: 
in a basis of sections on $D$, the meromorphic connection has matrix of connection 1-forms $\Om=-A(z)\ed z$, where the $n\times n$ matrix $A(z)$ equals
\[A(z)=\sum_{k=0}^\infty A_kz^{-k},\quad A_0\neq0.
\]Denote by $\C\{z^{-1}\}$ the ring of convergent power series in $z^{-1}$. The RHB Problem \ref{RHB1} is then equivalent to find a so-called {\it Birkhoff normal form}: does it exist a matrix $G\in GL(n,\C\{z^{-1}\})$ such that $B(z)=G^{-1}AG-G^{-1}\frac{d}{dz}G$ is of the form
\[B(z)=B_0+\frac{B_1}{z},\quad B_0,B_1\in M(n,\C)?
\]
\begin{rem}
The RHB Problem \ref{RHB1} is not always solvable. G.\,Birkhoff himself proved that the problem is solvable provided that the monodromy matrix of the differential system $\frac{d}{dz}Y(z)=A(z)Y$ is diagonalizable \cite{Bir13}, but he seemed to believe that the same would hold generally. This was disproved in 1959 by F.R.\,Gantmacher and P.\,Masani: they independently exhibited connections $\Om$ which cannot be put in Birkhoff normal form. See \cite{Gan59} and \cite{Mas59}. The counterexamples found by Gantmacher and Masani are of {\it reducible} nature, in the sense that they can be put in lower triangularly blocked form via an analytic transformation. This led to the following restricted problem: is the RHB Problem \ref{RHB1} solvable in the {\it irreducible} case? This question was answered positively, first for rank $n=2$ by W.B.\,Jurkat, D.A.\,Lutz, and A.\,Peyerimhoff \cite{JLP76}, then for $n=3$ by W.Balser \cite{Bal90}, and finally for any dimension by A.A.\,Bolibruch \cite{Bol94a,Bol94b}. We also refer to \cite[Ch.IV]{Sab07}, where the reader can find further {\it irreducibility} assumptions (named after J.\,Plemelj, A.A.\,Bolibruch, and V.\,Kostov) ensuring the solvability of the RHB Problem \ref{RHB1}. For the reducible case, the reader can find in \cite{BB97} sufficient conditions for the solvability of the problem.

\end{rem}
\begin{rem}The RHB Problem \ref{RHB1} admits several variants.
\begin{enumerate}
\item If one allows $B(z)$ to have a more general form $$B(z)=B_{-N}z^N+\dots+B_{-1}z+B_0+\frac{B_1}{z},\quad N\geq 1,\quad B_j\in M(n,\C),$$ then the problem always admits a positive solution. This was the original result proved by G.\,Birkhoff in \cite{Bir09}, see also \cite[\S 3.3]{Sib90}. 
\item If one allows {\it meromorphic} equivalences, i.e.\,\,gauge transformations $B(z)=G^{-1}AG-G^{-1}\frac{d}{dz}G$ with $G\in GL(n,\C\{z^{-1}\}[z])$, then the problem is known to be solvable in several cases. For example, if $n=2,3$, then the problem is always solvable as proved by W.B.\,Jurkat, D.A.\,Lutz, and A.\,Peyerimhoff \cite{JLP76}, and W.\,Balser \cite{Bal89}. For arbitrary ranks $n$, but under the assumption that $A_0$ has pairwise distinct eigenvalues, H.\,Turrittin showed that the problem always admits a positive solution \cite{Tur63} \cite[\S 3.10]{Sib90}.
\end{enumerate}
\end{rem}

\subsection{Families of Riemann--Hilbert--Birkhoff problems} Throughout the remaining part of the paper, $X$ will denote a connected complex manifold of dimension $d$. If $Z\subseteq X$ is a smooth analytic hypersurface, we denote by $\mathscr O_X(*Z)$ the sheaf of meromorphic functions on $X$ with poles on $Z$ at most. If $E$ is a holomorphic vector bundle on $X$, with sheaf of sections $\mathscr E$, we set $\mathscr E(*Z):=\mathscr E\otimes_{\mathscr O_X} \mathscr O_X(*Z)$.
\vskip2mm
In what follows, we want to consider families of RHB problems, parametrized by points of $X$.

\begin{defn}
Let $(E^o,\nabla^o)$ be a holomorphic vector bundle on a disc $D\subseteq \Pb^1$, centered at $z=\infty$, equipped with a meromorphic connection with a pole of order 2 at $z=\infty$. An {\it integrable deformation} $(\nabla,E,X,x_o)$ of $(E^o,\nabla^o)$ parametrized by $X$ is the datum of
\begin{itemize}
\item a vector bundle $E$ on $D\times X$,
\item a {\it flat} connection $\nabla$ on $E$ with a pole of order 2 along $\{\infty\}\times X$,
\item a point $x_o\in X$ at which $(E,\nabla)$ restricts to $(E^o,\nabla^o)$.
\end{itemize}
The integrable deformation is called {\it versal} if any other deformation with base space $X'$ is induced by the previous one via pull-back by a holomorphic map $\varphi\colon (X',x_o')\to (X,x_o)$. It is {\it universal} if the germ at $x_o'$ of the base-change $\varphi$ is unique. 
\end{defn}

Assume $(E^o,\nabla^o)$ to be extendable to a solution of the RHB Problem \ref{RHB1}: this means that (in a suitable basis of sections) the matrix of connections 1-forms of $\nabla^o$ takes the form 
\beq
\Om_o=-\left(A_o+\frac{B_o}{z}\right)\ed z.
\eeq 
Let $(\nabla,E,X,x_o)$ be an integrable deformation of $(E^o,\nabla^o)$. The next result shows that, for generic $x\in X$, the restriction $(E,\nabla)|_{D\times{x}}$ is extendable to a solution of the RHB Problem \ref{RHB1}, provided that $X$ is simply connected.

\begin{thm}[{\cite[Th.\,VI.2.1]{Sab07}\cite[Th.\,5.1(c)]{DH21}}]\label{thimp}Under the assumptions above, if $X$ is simply connected, then there exists 
\begin{itemize}
\item an analytic hypersurface $\Theta\subseteq X\setminus\{x_o\}$,
\item a \emph{unique} basis of sections of $\mathscr E(*(D\times \Theta))$,
\end{itemize}
with respect to which the the matrix of connection 1-forms of $\nabla$ takes the form
\beq\label{condef}
\Om=-\left(A(x)+\frac{B_o}{z}\right)\ed z-z\,C(x),\qquad x\in X\setminus \Theta,
\eeq
where
\begin{itemize}
\item the matrix $A(x)$ is a matrix of holomorphic functions on $X\setminus \Theta$, and meromorphic along $\Theta$, such that $A(x_o)=A_o$,
\item the matrix $C(x)$ is a matrix of holomorphic 1-forms on $X\setminus \Theta$, and meromorphic along $\Theta$, such that $C(x_o)=0$.\qed
\end{itemize}
\end{thm}
\begin{defn}
The matrix-valued holomorphic function $A\colon X\setminus \Theta\to M(n,\C)$ above is the {\it pole part} 
of the integrable deformation $\nabla$.

The matrix-valued 1-form $C\colon X\setminus \Theta\to M(n,\C)\otimes \Om^1_X$ above is the {\it deformation part} of the integrable deformation $\nabla$.
\end{defn}
The integrability condition for $\nabla$ translates into the system of equations 
\beq
\label{inteq1}
\ed C=0,\quad C\wedge C=0,\quad [A,C]=0,\quad \ed A=C+[C,B_o].
\eeq
If $\bm x=(x^1,\dots, x^d)$ are local holomorphic coordinates on $X$, and if $C=\sum_{i=1}^dC_i(\bm x)\ed x^i$, these equations take the form
\beq
\frac{\der C_i}{\der x^j}=\frac{\der C_j}{\der x^i},\quad [C_i,C_j]=0,\quad [A,C_i]=0,\quad \frac{\der A}{\der x^i}=C_i+[C_i,B_o],
\eeq for $i,j=1,\dots,d$.

\subsection{Universal integrable deformations: Malgrange's and Jimbo--Miwa--Ueno's theorems}\label{secmal}
Let $(E^o,\nabla^o)$ be a solution of a RHB problem \ref{RHB1}, i.e. a trivial vector bundle (of rank $n$) on $\Pb^1$ with meromorphic connection with matrix (in a suitable basis of sections) of the form
\beq\label{con1}\Om_o=-\left(A_o+\frac{B_o}{z}\right)\ed z.
\eeq
\vskip2mm
\noindent{\bf Question 1:} Under which conditions on $(A_o,B_o)$ does there exist a universal integrable deformation of the connection \eqref{con1}?
\vskip2mm
A first positive result, due to B.\,Malgrange, requires that the pole part $A_o$ is an element of the regular stratum $\Mreg$ of Section \ref{secmreg}.

\begin{thm}[{\cite{Mal83a,Mal86}}]\label{thmal1}Assume that the matrix $A_o$ is regular. The connection $\nabla^o$  with matrix \eqref{con1} has a germ of universal deformation.
\end{thm}
The reader can find the proof in Appendix \ref{appmal}.
\vskip2mm
This result can be made more explicit, 
under the further {\it semisimplicity} assumption on $A_o$. Assume that $A_o\in\Mreg\cap\Mdiag=\mc M_{\msl1;1;\dots;1\msr}$, and let $P\in GL(n,\C)$ such that 
\[P^{-1}A_oP=\La_o={\rm diag}(u^1_o,\dots, u^n_o),\quad u_o^i\neq u_o^j,\text{ for }i\neq j.
\]Set
\beq
\label{con1.2}
\widehat\Om_o:=P^{-1}\Om_o P=-\left(\La_o+\frac{\mc B_o}{z}\right)\ed z,\quad \mc B_o:=P^{-1}B_oP.
\eeq 

For $\bm u\in\C^n$, denote $\La(\bm u):=\diag(u^1,\dots,u^n)$, so that $\La(\bm u_o)=\La_o$. Given a matrix $M$ denote by $M'$ its diagonal part, and by $M''$ its off-diagonal part.
\begin{thm}[\cite{Mal83b,Mal86}]\label{thmal2}
Under the assumptions above, there exists 
a sufficiently small polydisc $\mathbb D=\mathbb D(\bm u_o)\subseteq \C^n$ with center at $\bm u_o$, 
and a holomorphic off-diagonal matrix $\Gm\colon \mathbb D\to M(n,\C)$, $\Gm(\bm u)=\Gm''(\bm u)$, 
such that:
\begin{enumerate}
\item  the matrix of 1-forms $\widehat\Om$ on $\C^*\times \mathbb D$, defined by
\beq\label{uic}
\widehat\Om(z,\bm u):=-\ed\left(z\La(\bm u)\right)-\left([\Gm(\bm u),\La(\bm u)]+\mc B'_o\right)\frac{\ed z}{z}-[\Gm(\bm u),\ed \La(\bm u)],
\eeq
defines an integrable connection $\nabla$ on the trivial bundle $\underline{\C^n}\to \C^*\times \mathbb D$;
\item the $\ed z$-component of $\widehat\Om$ restricts to $\widehat\Om_o$ at $\bm u_o$, i.e. 
\beq\label{eqresom}
\widehat\Om(z,\bm u_o)=\widehat\Om_o+\om,\quad \om\in M(n,\Om^1_{\mathbb D});
\eeq
\item $\nabla$ is formally equivalent at $z=\infty$ to the matrix connection
\beq\label{con2}
-\ed\left(z\La(\bm u)\right)-\mc B'_o\frac{\ed z}{z},
\eeq
that is there exists a $z^{-1}$-formal base change $\Phi(z,\bm u)=\sum_{k=0}^\infty\Phi_k(\bm u)z^{-k}$, with $\Phi_k\colon\mathbb D\to M(n,\C)$ holomorphic and $\Phi_0(\mathbb D)\subseteq GL(n,\C)$, such that
\[
\Phi^{-1}\widehat\Om\Phi+\Phi^{-1}\ed \Phi=-\ed\left(z\La(\bm u)\right)-\mc B'_o\frac{\ed z}{z};
\]
\item $\nabla$ defines a \emph{universal} integrable deformation of its restriction at any point $\bm u\in\mathbb D$.
\end{enumerate}
Moreover, the matrix $\Gm$ is \emph{uniquely} determined by these conditions. \qed
\end{thm}
\begin{rem}The integrability condition of $\nabla$ is equivalent to the following equations 
\[
\ed [\Gm,\ed\La]=[\Gm,\ed\La]\wedge [\Gm,\ed\La],\qquad
\ed[\Gm,\La]=\left[[\Gm,\ed\La],\mc B'_o+[\Gm,\La]\right],
\]called {\it Darboux--Egoroff equations.} In local coordinates $\bm u$, they read
\begin{align}
\label{DE01}
\der_k\Gm_{ij}&=\Gm_{ik}\Gm_{kj},\quad k\neq i,j,\\
\label{DE02}
(u^j-u^i)\der_i\Gm_{ij}&=\sum_{k\neq i,j}(u^k-u^j)\Gm_{ik}\Gm_{kj}-(b_j-b_i-1)\Gm_{ij},\\
\label{DE03}
(u^i-u^j)\der_{j}\Gm_{ij}&=\sum_{k\neq i,j}(u^k-u^i)\Gm_{ik}\Gm_{kj}-(b_j-b_i-1)\Gm_{ij},
\end{align}where $\Gm=(\Gm_{ij})_{i,j=1}^n$, and $\mc B_o'={\rm diag}(b_1,\dots, b_n)$.
\end{rem}
The statement of Theorem \ref{thmal2} ca be further refined to a global one.
Let $\Dl$ be the union of {big diagonal} hyperplanes in $\C^n$, defined by the equations
\[\Dl:=\bigcup_{i<j}\{\bm u\in\C^n\colon u^i=u^j\},
\]let $X_n$ be the complement $\C^n\setminus \Dl$, with base point $\bm u_o:=(u_o^1,\dots, u_o^n)$. Denote by $\pi\colon (\widetilde X_n,\tilde{ \bm u}_o)\to (X_n,\bm u_o)$ the universal cover of $X_n$, equipped with fixed base points $\tilde {\bm u}_o$ and $\bm u_o$, respectively. The space $X_n$ is identified with the space of diagonal regular $n\times n$ matrices.

\begin{thm}[\cite{JMU81,Mal83b}]\label{thmal}
There exists on $\Pb^1\times \widetilde X_n$ a vector bundle $E$, equipped with a meromorphic connection $\nabla$, such that
\begin{enumerate}
\item the coefficients of $\nabla$ have poles along the hypersurface $\Theta\subseteq \widetilde X_n$ of points $\tilde{\bm u}\in\widetilde X_n$ such that $E|_{\Pb^1\times\{\tilde{\bm u}\}}$ is not trivial;
\item  $\nabla$ is flat, with a pole of Poincar\'e rank 1 along $\{\infty\}\times \widetilde X_n$, and a logarithmic pole along $\{0\}\times \widetilde X_n$;
\item $(E,\nabla)$ restricts to $(E^o,\nabla^o)$ at $\tilde{\bm u}_o$;
\item for any $\tilde{\bm u}\in\widetilde X_n$, the eigenvalues of the pole part of $\nabla$ at the point $(\infty,\tilde{\bm u})$ equal (up to permutation) the $n$-tuple $\pi(\tilde{\bm u})$.
\end{enumerate}
Moreover, for any $\tilde{\bm u}\in\widetilde X_n\setminus\Theta$, the bundle with meromorphic connection $(E,\nabla)$ induces a universal deformation of its restriction $(E,\nabla)|_{\Pb^1\times\{\tilde{\bm u}\}}$. \qed
\end{thm}

\subsection{Integrable deformations of degenerate Birkhoff normal forms: Sabbah's theorem} Malgrange's and Jimbo--Miwa--Ueno's Theorems \ref{thmal1}, \ref{thmal2}, \ref{thmal} provide an answer to Question 1, in the case the pole part $A_o\in\Mreg$.  In a sense, these results are the best possible: if $A_o\notin \Mreg$, then in general the connection \eqref{con1} does not admit versal deformations, see the example in Appendix \ref{appnover}.
\vskip2mm
Consider now the stratum $\Mdiag$ of diagonalizable matrices, that is the conjugate stratum of $\Mreg$ in the sense of Section \ref{secmreg}. 

Let us assume that $A_o\in\Mdiag$: 
in the notations of the previous section, assume $\bm u_o\in\Dl$. Define the partition $\{1,\dots,n\}=\coprod_{r\in R}I_r$ such that for any $r\in R$ we have
\[\{i,j\}\subseteq I_r\quad\text{if and only if}\quad u_o^i=u_o^j.
\]  In \cite{Sab21}, C.\,Sabbah addressed the following problem. 
\vskip2mm
\noindent{\bf Question 2: }Is it possible to find an integrable deformation of the form \eqref{uic} of the Birkhoff normal form \eqref{con1.2} with $z^{-1}$-formal normal form \eqref{con2}? 
\vskip2mm
Remarkably, in \cite[Section 4]{Sab21} it is shown that the answer is positive, under (sharp) sufficient conditions on the coefficient $\mc B_o$ of the normal form \eqref{con1.2}.
\vskip2mm
\noindent{\bf Property PNR\footnote{PNR stands for ``partial non-resonance''.}: }There exists a matrix $P\in GL(n,\C)$ diagonalizing $A_o$, i.e. $P^{-1}A_oP=\La_o={\rm diag}(u_o^1,\dots, u_o^n)$, and such that the matrix $\mc B_o:=P^{-1}B_oP$ has the following properties:
\begin{itemize}
\item[$(\star)$] $\mc B''_o\in {\rm Im\ ad}(\La(\bm u_o))$.
\item[$(\star\star)$] $\mc B'_o$ is \emph{partially non-resonant}, i.e. 
\[\forall\ r\in R,\quad \forall\ i,j\in I_r,\quad (\mc B'_o)_{ii}-(\mc B'_o)_{jj}\notin \Z\setminus\{0\}.\]
\end{itemize}
\begin{thm}[{\cite[Th. 4.9]{Sab21}}]\label{tsab}
Assume that Property PNR above holds true. 
Let $\bm u_o\in\Dl$, and $\mc V$ a neighborhood of $\bm u_o$ in $\C^n$. 
If $\mc V$ is sufficiently small, there exists a holomorphic hypersurface $\Theta$ in $\mc V\setminus\{\bm u_o\}$ and a holomorphic off-diagonal matrix $\Gm''(\bm u)$ on $\mc V\setminus\Theta$, such that
\begin{enumerate} 
\item the 1-forms matrix \eqref{uic} defines a meromorphic connection $\nabla$ on the trivial vector bundle on $\mathbb P^1\times (\mc V\setminus\Theta)$; 
\item $\nabla$ restricts to the connection \eqref{con1.2} at $\bm u_o$;
\item $\nabla$ is formally equivalent at $z=\infty$ to the connection defined by the matrix of 1-forms \eqref{con2}. 
\end{enumerate}The matrix $\Gm''$ is \emph{uniquely} determined by these conditions. \qed
\end{thm}
\begin{rem}
Property PNR-$(\star)$ is equivalent to $(\mc B_o)_{ij}=0$ whenever $i,j\in I_r$ for some $r$. It is a necessary condition for the statement of Theorem \ref{tsab}: by restriction of \eqref{uic} at $\bm u=\bm u_o$, we obtain $\mc B_o''=[\Gm(\bm u_o),\La_o]$. 
\end{rem}

\subsection{Integrable deformations of $\rm d$/$\rm dv$/$\rm fs\,$-type} Let $\nabla^o$ be a connection on a trivial vector bundle $E^o\to\Pb^1$with matrix of connection 1-forms \eqref{con1}. Consider an integrable deformation $(\nabla,E,X,x_o)$ of $\nabla^o$, parametrized by a complex manifold $X$, with matrix of connection 1-forms $\Om$ as in equation \eqref{condef}. In particular let 
\[A\colon X\setminus \Theta\to M(n,\C),\qquad C\colon X\setminus\Theta\to M(n,\C)\otimes \Om^1_X,
\]be the pole and deformation parts, respectively.

\begin{defn}
The deformation $(\nabla,E,X,x_o)$ is said to be 
\begin{itemize}
\item of {\it diagonal type} (for short, d{\it-type}) if the pole part $A$ and the deformation part $C$ of $\Om$ are locally holomorphically diagonalizable matrices at $x_o$;
\item of {\it generic diagonal type} (for short, {\it generic $\rm d$-type}) if it is of $\rm d$-type, and if $f_1,\dots, f_n$ are the holomorphic eigenvalues of $A$, then we have $\ed_{x_o}f_i \neq \ed_{x_o}f_j$, for any $i\neq j$. 
\end{itemize}
A germ of integrable deformation will be said of (generic) d-type if at least one (and hence any) of its representative is of (generic) d-type.
We denote by $\frak I_{\rm d}(\nabla^o)$ (resp. $\frak I_{\rm d}^{\rm gen}(\nabla^o)$) the classes of germs of integrable deformations of $\nabla^o$ which are of $\rm d$-type (resp. generic $\rm d$-type).
\end{defn}

\begin{rem}
The genericity condition implies that no 1-form $\ed(f_i-f_j)$ is vanishing in a neighborhood of $x_o$. In particular, no pair of eigenvalues $f_i$ and $f_j$, with $i\neq j$, are identically equal.
\end{rem}

\begin{thm}\label{thforcan1}
Let $(\nabla,E,X,x_o)$ be an integrable deformation of $\nabla^o$ of \emph{d}-type. Let $U\subseteq X\setminus\Theta$ a neighborhood of $x_o$ such that $A|_U,\,C|_U$ are holomorphically diagonalizable.
\begin{enumerate}
\item If $\Dl_0\colon U\to M(n,\C)$, with $\Dl_0''=0$, is the holomorphic diagonal form of $A|_U$, then $\ed \Dl_0\in\Om^1_U\otimes M(n,\C)$ is the diagonal form of $C|_U$.
\item There exists a base of holomorphic sections of $E|_{\Pb^1\times U}$ with respect to which $\nabla$ has the following matrix $\widetilde{\Om}$ of connection 1-forms
\beq\label{forcan2}
\widetilde\Om(z,x)=-\left(\Dl_0(x)+\frac{1}{z}\mc B(x)\right)\ed z-z\,\ed \Dl_0(x)+\varpi(x),
\eeq where $\mc B\colon U\to M(n,\C)$ is holomorphic, and $\varpi\in\Om^1_U\otimes M(n,\C)$. 
\end{enumerate}
\end{thm}

\begin{rem}\label{remintd}
The integrability conditions for $\nabla$, in terms of $(\Dl_0,\mc B,\varpi)$ as in \eqref{forcan2}, read
\[[\ed \Dl_0,\mc B]+[\Dl_0,\varpi]=0,\qquad \ed \mc B=[\mc B,\varpi],\qquad \ed \Dl_0\wedge\varpi+\varpi\wedge\ed \Dl_0=0,\qquad \ed\varpi+\varpi\wedge\varpi=0.
\]Moreover, notice that $\mc B$ is holomorphically similar to the constant matrix $B_o$, and hence holomorphically Jordanizable.
\end{rem}

\proof Consider the matrix $\Om(z,x)=-\left(A(x)+\frac{1}{z}B_o\right)\ed z-z C(x)$ defining $\nabla$. The matrices $A$ and $C$ are, by assumption, locally holomorphically diagonalizable at $x_o$. Moreover we have $[A,C]=0$, by integrability of $\nabla$. Consequently, there exists a holomorphic map $G\colon U\to GL(n,\C)$ which {\it simultaneously} diagonalizes $A(x)$ and $C(x)$, i.e.
\[G(x)^{-1}A(x)G(x)=\Dl_0(x),\qquad G(x)^{-1}C(x)G(x)=\Dl_1(x). 
\]
Set  $\mc B:=G^{-1}B_oG$, and $\varpi=G^{-1}\ed G$. We have
\[G^{-1}\Om G+G^{-1}\ed G=-\left(\Dl_0(x)+\frac{1}{z}\mc B\right)\ed z-z\Dl_1(x)+\varpi,
\]and the integrability conditions read
\[\ed \Dl_0=\Dl_1+[\Dl_1,\mc B]+[\Dl_0,\varpi],\qquad \ed \mc B=[\mc B,\varpi],\qquad [\Dl_0,\Dl_1]=0,
\]
\[\ed \Dl_1+\Dl_1\wedge\varpi+\varpi\wedge\Dl_1=0,\qquad \Dl_1\wedge\Dl_1=0,\qquad \ed \varpi+\varpi\wedge\varpi=0.
\]
From the first equation we deduce $\Dl_1=\ed \Dl_0$. 
\endproof 

\begin{defn}
Let $(\nabla,E,X,x_o)$ be an integrable deformation of d-type. We say that $(\nabla,E,X,x_o)$ is of {\it diagonal-vanishing type} (for short, dv-{\it type}) if there exist
\begin{itemize}
\item a neighborhood $U\subseteq X\setminus\Theta$ of $x_o$,
\item a holomorphic off-diagonal matrix $\mc L\colon U\to M(n,\C)$, $\mc L''=\mc L$,
\item a basis of holomorphic sections of $E|_{\Pb^1\times U}$,
\end{itemize}with respect to which $\nabla$ has matrix of connection 1-forms as in \eqref{forcan2} with
\beq
\label{veq}
\mc B''=[\mc L,\Dl_0],\qquad \varpi''=[\ed \Dl_0,\mc L].
\eeq
We say that a germ of integrable deformation is of dv-type if at least one (and hence any) of its representative is of dv-type. We denote by $\frak I_{\rm dv}(\nabla^o)$ the class of germs of dv-type integrable deformations of $\nabla^o$.
\end{defn}

\begin{thm}\label{thDL}
Let $(\nabla,E,X,x_o)$ be an integrable deformation of $\nabla^o$ of ${\rm dv}$-type. There exist a neighborhood $U$ of $x_o$, and a basis of sections of $E|_{\Pb^1\times U}$ with respect to which $\nabla$ has the following matrix of connection 1-forms 
\beq\label{forcan3}
\widehat\Om(z,x)=-\left(\Dl_0(x)+\frac{1}{z}\frak B(x)\right)\ed z-z\,\ed \Dl_0(x)+\om(x),
\eeq
where
\beq\label{forcan3.2}
\frak B'=\frak B'_o=\text{const.},\qquad \frak B''=[L,\Dl_0],\qquad \om'=0,\qquad \om''=[\ed \Dl_0,L],
\eeq for some holomorphic off-diagonal matrix $L\colon U\to M(n,\C)$.
\end{thm}
\proof Let $U$ be as in Theorem \ref{thforcan1}. 
Consider a matrix $\widetilde \Om$ defining $\nabla$ as in \eqref{forcan2}. By splitting $\varpi=\varpi'+\varpi''$, we have
\[\ed \varpi'=-(\varpi''\wedge\varpi'')'=-\left([\ed \Dl_0,\mc L]\wedge [\ed \Dl_0,\mc L]\right)'=0.
\]Since $\varpi'$ is closed, locally there exists an invertible diagonal matrix $H\colon U\to GL(n,\C)$ such that $\varpi'=-H^{-1}\ed H$. The matrix $\widehat\Om=H^{-1}\widetilde\Om H+H^{-1}\ed H$ is as in \eqref{forcan3}, with
\[\frak B=H^{-1}\mc B H,\qquad \om=H^{-1}\varpi H+H^{-1}\ed H.
\]The last three equations of \eqref{forcan3.2} are automatically satisfied, with $L=H^{-1}\mc L H$. The first equation follows from the integrability condition $\ed \frak B=[\frak B,\om]$: we have
\[
\pushQED{\qed}
\ed\frak B'=[\frak B'',\om'']'=[[L,\Dl_0],[\ed \Dl_0,L]]'=0.\qedhere
\popQED
\]

\begin{rem}
In terms of the matrices $(\Dl_0,\frak B'_o, L)$, the integrability condition of the connection \eqref{forcan3} reads
\begin{align}
\label{DEeq1}
\ed [L,\ed\Dl_0]&=[L,\ed\Dl_0]\wedge [L,\ed\Dl_0],\\
\label{DEeq2}
\ed[L,\Dl_0]&=\left[[L,\ed\Dl_0],\frak B'_o+[L,\Dl_0]\right].
\end{align}
We call these equations the {\it generalized Darboux--Egoroff equations}.
\end{rem}

\begin{defn}
Let $(\nabla,E,X,x_o)$ be an integrable deformation of $\nabla^o$ with matrix $\Om$ as in \eqref{condef}. We say that $(\nabla,E,X,x_o)$ is of {\it formally simplifiable type} (for short, fs\,-{\it type}) if there exist
\begin{itemize}
\item a neighborhood $U\subseteq X\setminus\Theta$ of $x_o$, 
\item a sequence of holomorphic maps $\Phi_k\colon U\to M(n,\C)$, with $k\geq 0$ and $\Phi_0(U)\subseteq GL(n,\C)$, 
\item a holomorphic diagonal map $\Dl_0\colon U\to M(n,\C)$, 
\item a constant diagonal matrix $\frak B_o'\in M(n,\C)$,
\end{itemize}such that
\[\Phi^{-1}\Om\Phi+\Phi^{-1}\ed \Phi=-\ed \left(z\Dl_0(x)\right)-\frak B_o'\frac{\ed z}{z},\qquad \Phi:=\sum_{k\geq 0}\Phi_kz^{-k}.
\]We say that a germ of integrable deformation is of fs\,-type if at least one (and hence any) of its representative is of fs\,-type. We denote by $\frak I_{\rm fs}(\nabla^o)$ the class of germs of fs\,-type integrable deformations of $\nabla^o$.
\end{defn}

\begin{rem}\label{remimp}
The functions $\Phi_k$ in the definition of integrable deformation of fs\,-type satisfy the following equations:
\begin{align*}
&A\Phi_0=\Phi_0\Dl_0,\qquad C\Phi_0=\Phi_0\ed \Dl_0,\\ 
&A\Phi_{k+1}+B_o\Phi_k+k\Phi_k=\Phi_{k+1}\Dl_0+\Phi_k\frak B_o',\quad &k\geq 0,\\
&\ed\Phi_k=C\Phi_{k+1}-\Phi_{k+1}\ed \Dl_0,\quad &k\geq0.
\end{align*}
\end{rem}
\vskip2mm
Denote by $\frak I(\nabla^o)$ the set of all germs of integrable deformations of a connection $\nabla^o$.
\begin{thm}\label{mth1} $\qquad$
\begin{enumerate}
\item The classes $\Id,\Idv,$ and $\Ifs$ are closed by arbitrary base change of the deformation parameter spaces $\phi\colon (X,x_o)\to (X',x_o')$.
\item For any connection $\nabla^o$ of the form \eqref{con1}, we have
\[\xymatrix @R=.01pc@C=0.1pc{
\Ifs_{\dsubp}&\\
&\!\!\!\!\!\!\Idv \subseteq \Id \subseteq \frak I(\nabla^o)\\
\Igend^{\dsub}&
}
\]
\item If the pole part $A_o$ of $\nabla^o$ is in $\Mreg\cap\Mdiag$, then we have
\[\Ifs = \Idv = \Id = \frak I(\nabla^o).
\]
\end{enumerate}
\end{thm}
\proof
Point (1) is obvious, by definition of the classes $\Id,\Idv$, and $\Ifs$.
\vskip2mm
The only nontrivial inclusion of point (2) is $(\Ifs\cup\Igend)\subseteq \Idv$. Let $\nabla\in\Ifs$ be defined by the matrix $\Om$ as in \eqref{condef}. Let $(\Phi_k)_{k\geq 0}$ be the matrix-valued functions such that $\Phi:=\sum_{k\geq0}\Phi_kz^{-k}$ satisfies $\Phi^{-1}\Om\Phi+\Phi^{-1}\ed \Phi=-\ed \left(z\Dl_0(x)\right)-\frak B_o'\frac{\ed z}{z}$. By invoking the equations of Remark \ref{remimp}, it is easy to see that
\[\Phi_0^{-1}\Om\Phi_0+\Phi^{-1}_0\ed \Phi_0=-\left(\Dl_0+\frac{1}{z}(\frak B_o'+[\Phi_0^{-1}\Phi_1,\Dl_0])\right)\ed z-z\ed \Dl_0-[\Phi_0^{-1}\Phi_1,\ed\Dl_0].
\]This proves that $\Ifs\subseteq\Idv$. 

Let $\nabla\in\Igend$.  By Theorem \ref{thforcan1}, it can be defined by a matrix $\widetilde\Om$ of connection 1-forms as in \eqref{forcan2} for a suitable pair $(\mc B,\varpi)$ of matrices. By integrability, we have $\ed\Dl_0\wedge\varpi+\varpi\wedge\ed\Dl_0=0$.  The equation for the $(i,j)$ entry, with $i\neq j$, reads $\varpi_{ij}\wedge(\ed f_i-\ed f_j)=0$. Assume temporarily that the deformation is parametrized by a complex manifold $X$ with $\dim_\C X\geq 2$. Since $\ed f_i-\ed f_j\not\equiv 0$ in a neighborhood $U$ of $x_o$, we deduce the existence\footnote{If $\dim_\C X=1$, the vanishing condition $\varpi_{ij}\wedge(\ed f_i-\ed f_j)=0$ does not imply that $\varpi_{ij}$ is a multiple of $(\ed f_i-\ed f_j)$. The vanishing condition is indeed satisfied by any arbitrary holomorphic 1-form $\varpi_{ij}$. } of a holomorphic function $\mc L_{ij}\colon U\to \C$ such that $\varpi_{ij}=\mc L_{ij}(\ed f_i-\ed f_j)$. This is a special case of de Rham's division lemma, see \cite{dR54}\cite[Lemma 3.1]{NY04}. By integrability, we also have $[\ed\Dl_0,\mc B]+[\Dl_0,\varpi]=0$. The equation for the $(i,j)$ entry reads $(\ed f_i-\ed f_j)\mc B_{ij}+(f_i-f_j)\varpi_{ij}=0$.  Hence, we necessarily have $\mc B_{ij}+(f_i-f_j)\mc L_{ij}=0$. This shows that, if $\dim_\C X\geq2$, $\nabla$ is of dv-type, since $\mc B''=[\mc L,\Dl_0]$ and $\varpi''=[\ed \Dl_0,\mc L]$. 

The case $\dim_\C X=1$ can be reduced to the case $\dim_\C X>1$. If $(\nabla,E,X,x_o)$ is a (germ of) integrable deformation of (generic) d-type, we can construct a new (germ of) integrable deformation $(\nabla',E',X\times\C,(x_o,0))$ as follows. Set $E':={\rm pr}^*E$, where ${\rm pr}\colon X\times\C\to X$. If $\Om(z,x)$ is the matrix of $\nabla$ as in \eqref{condef}, define $\nabla'$ by the matrix of connection 1-forms 
\[\Om'(z,x,s)=\Om(z,x)+\ed(zs\cdot{\rm Id}_{n}),\qquad (x,s)\in X\times \C,\quad n={\rm rk}(E).
\]We have $\nabla=\iota^*\nabla'$, where $\iota\colon X\to X\times\C$ is the canonical embedding $\iota(x)=(x,0)$. 
The connection $\nabla'$ is of generic d-type, so that $\nabla'\in\Idv$. The result for $\nabla$ follows from point (1). This proves point (2).
\vskip2mm
Let us prove point (3). Consider an arbitrary integrable deformation defined by a matrix $\Om$ as in \eqref{condef}. The set $\Mreg\cap\Mdiag$ is an open dense subset of $M(n,\C)$. Hence, 
if $A_o\in\Mreg\cap\Mdiag$, there exists a neighborhood $U\subseteq X$ of $x_o$ such that $A(x)\in\Mreg\cap\Mdiag$ for any $x\in U$. Consequently, $A(x)$ is locally holomorphically diagonalizable by Theorem \ref{thjb2}. By integrability, the matrices $A(x)$ and $C(x)$ commute for any $x\in U$. By regularity of $A(x)$, it follows that $C(x)$ is a polynomial expression of $A(x)$. Thus, $C(x)$ is locally holomorphically diagonalizable. This shows that $\frak I(\nabla^o)=\Id$. 

Consider an integrable deformation $\nabla$ of d-type, defined by a matrix $\widetilde \Om$ as in \eqref{forcan2} for a suitable pair $(\mc B,\varpi)$. By integrability, we have $[\ed\Dl_0,\mc B]+[\Dl_0,\varpi]=0$. The equation for the $(i,j)$ entry, with $i\neq j$, reads $(\ed f_i-\ed f_j)\mc B_{ij}+(f_i-f_j)\varpi_{ij}=0$. By the argument above, we also have $f_i(x)\neq f_j(x)$ for any $i\neq j$ and any $x\in U$. Hence, if we set $\mc L_{ij}:=\frac{\mc B_{ij}}{f_j-f_i}$, we obtain a holomorphic off-diagonal matrix $\mc L\colon U\to M(n,\C)$ such that $\mc B''=[\mc L,\Dl_0]$ and $\varpi''=[\ed\Dl_0,\mc L]$. This shows that $\Id=\Idv$.

It thus remains to show that $\Idv\subseteq\Ifs$. Let $\nabla\in\Idv$ to be defined by the matrix $\widehat\Om$ as in \eqref{forcan3}. As before, we can take $U$ sufficiently small so that $f_i(x)\neq f_j(x)$ for any $i\neq j$ and any $x\in U$. We claim that there exists a unique sequence of holomorphic functions $F_k\colon U\to M(n,\C)$, with $k\geq 1$, such that
\[F^{-1}\widehat\Om F+F^{-1}\ed F=-\ed \left(z\Dl_0(x)\right)-\frak B_o'\frac{\ed z}{z},\]where\[ F(z,x):={\rm Id}_n+\sum_{k=1}^\infty F_k(x)z^{-k},
\qquad \frak B'_o=\frak B(x_o)'.
\]This holds if and only if the following identities are satisfied
\begin{align}
\label{ie1}&[\Dl_0,F_{k+1}]+\frak BF_k-F_k\frak B'_o+kF_k=0,\qquad &&k\geq 0,	\quad F_0={\rm Id}_n,\\
\label{ie2}&\ed F_{k}=[F_1,\ed\Dl_0]F_k+[\ed \Dl_0,F_{k+1}],\qquad &&k\geq 0.
\end{align}
The functions $F_k$'s can be found by iteratively solving equations \eqref{ie1}.
For $i\neq j$, we find 
\[(F_{k+1})_{ij}=\frac{1}{f_j-f_i}\left(\frak BF_k-F_k\frak B'_o+kF_k\right)_{ij},\qquad k\geq 0
\]and the diagonal elements are given by
\[(F_{k+1})_{ii}=-\frac{1}{k+1}\sum_{\ell\neq i}\frak B_{i\ell}(F_{k+1})_{\ell i},\qquad k\geq 0.
\]The resulting functions $F_k$'s are holomorphic on $U$. In particular, notice that $F_1''=L$. Let us prove that equations \eqref{ie2} are automatically satisfied. The proof is standard, but for completeness we outline it below. Fix the sector $\mc V$ in the universal cover of $\C^*$ defined by $\mc V=\{z\in\widetilde{\C^*}\colon |\arg z|<\frac{\pi}{2}\}$. By a theorem of Y.\,Sibuya \cite[Main Th.]{Sib62} \cite{HS66} \cite{Was65,BJL79}, the differential system of equations 
\beq
\label{ode1}
\frac{d}{dz}Y=\left(\Dl_0(x)+\frac{1}{z}\frak B(x)\right)Y
\eeq
admits a unique fundamental system of solutions $Y(z,x)$ such that 
\beq
\label{asym}
Y(z,x)z^{-\frak B_o'}e^{-\Dl_0(x)}\sim F(z,x),\quad |z|\to\infty,\,z\in\mc V,\quad \text{uniformly in }x.
\eeq
If we prove that $Y$ satisfies the system of equations
\beq\label{syst2}
\ed' Y=(z\,\ed'\Dl_0(x)+[F_1,\ed'\Dl_0])\,Y, \quad \ed':=\sum_{i=1}^d\frac{\der}{\der x^i}\ed x^i,\quad x^i\text{'s coordinates on $X$,}
\eeq
then equations \eqref{ie2} immediately follow from the asymptotic expansion \eqref{asym}. For any $i=1,\dots, d$, set \[W_i(z,x):=\der_iY(z,x)-(z\der_i\Dl_0(x)+[F_1,\der_i\Dl_0])\,Y(z,x).\]
We claim that each function $W_i$ is a solution of \eqref{ode1}: by a simple computation, and by invoking the first two identities of Remark \ref{remintd}, we find that
\[\der_zW_i-\left(\Dl_0+\frac{1}{z}\frak B\right)W_i=\der_i\left(\der_zY-\left(\Dl_0+\frac{1}{z}\frak B\right)Y\right)=0.
\]Hence, for any $x\in X$, there exists a matrix $C_i(x)$ such that $W_i(z,x)=Y(z,x)C_i(x)$. On the one hand, we have
\[W_i(z,x)\sim\left(\der_iF+z[F,\der_i\Dl_0]+[F_1,\der_i\Dl_0]F\right)z^{\frak B_o'}e^{z\Dl_0},\quad |z|\to\infty,\,z\in\mc V,\quad \text{ uniformly in }x.
\]On the other hand, we also have 
\[W_i(z,x)\sim F(z,x)z^{\frak B_o'}e^{z\Dl_0}C_i(x),\quad |z|\to\infty,\,z\in\mc V,\quad \text{ uniformly in }x.
\]Consequently, we deduce that
\beq
\label{asym2}
z^{\frak B_o'}e^{z\Dl_0}C_i(x)e^{-z\Dl_0}z^{-\frak B_o'}=\text{formal power series in }\frac{1}{z}.
\eeq
For a fixed $x\in X$, and for $j\neq k$, the sector $\mc V$ contains rays $\ell$ of points $z$ along which ${\rm Re}(z(f_j(x)-f_k(x)))>0$. We deduce that the $(j,k)$-entry of $C_i(x)$ vanishes, otherwise we would have a divergence, for $|z|\to \infty$ along the rays $\ell$, on the l.h.s. of \eqref{asym2}. So the matrix $C_i(x)$ is diagonal, and
we have
\begin{align*}C_i(x)&=z^{\frak B_o'}e^{z\Dl_0}C_i(x)e^{-z\Dl_0}z^{-\frak B_o'}\\
&=F(z,x)^{-1}\left(\der_iF(z,x)+z[F(z,x),\der_i\Dl_0(x)]-[F_1(x),\der_i\Dl_0(x)]F(z,x)\right)\\
&=z(\der_i\Dl_0-\der_i\Dl_0)+\left(F_1\,\der_i\Dl_0-\der_i\Dl_0\,F_1-[F_1,\der_i\Dl_0]\right)+O\left(\frac{1}{z}\right)=O\left(\frac{1}{z}\right).
\end{align*}
This shows that $C_i(x)=0$ for any $i=1,\dots, d$. Hence, \eqref{syst2} hold true. This completes the proof. \endproof

\begin{rem}\label{rkstrict}
In point (2) of Theorem \ref{mth1}, we expect the inclusions of the chain $\Ifs\subseteq\Idv\subseteq\Id$ to be strict. It would be interesting to show this via explicit examples.
\end{rem}

The following result further clarifies the relation between the classes $\Ifs$ and $\Igend$.

\begin{thm}\label{mth2}Let $\nabla^o$ be a connection with matrix \eqref{con1}.
Assume one of the following assumptions hold:
\begin{itemize}
\item The pole part $A_o$ is an element of $\Mreg\cap\Mdiag$;
\item The pole part $A_o$ is an element of $\Mdiag$ and the Property PNR holds.
\end{itemize}
Then we have $\emptyset\neq\Igend\subsetneq\Ifs$.
\end{thm}

\proof
Let us first assume that $A_o\in\Mdiag\cap\Mreg$. We have $\Igend\subseteq\Ifs$, by point (3) of Theorem \ref{mth1}. Consider the universal integrable deformation $(\nabla^{\rm JMUM},\underline{\C^n},\mathbb D(\bm u_o),\bm u_o)$ of $\nabla^o$, whose existence is established by Malgrange's and Jimbo--Miwa--Ueno's Theorems \ref{thmal1} and \ref{thmal}. Its germ is an element of $\Igend$. Moreover, if $\Om_o(z)$ is the matrix of connection 1-forms defining $\nabla^o$ as in \eqref{con1}, then the matrix $\Om(z,s):=\Om_o(z)+\ed(zs\cdot {\rm Id}_n)$, with $s\in\C$, defines an element of $\Ifs\setminus\Igend$. This proves the statement in the case $A_o$ is in $\Mreg\cap\Mdiag$.
\vskip2mm
Let us now assume that $A\in\Mdiag$, and that the Property PNR holds. Consider the integrable deformation $(\nabla^{\rm Sab},\underline{\C^n},\mathbb D(\bm u_o),\bm u_o)$ of $\nabla^o$, whose existence is established by Sabbah's Theorem \ref{tsab}. It defines an element of $\Ifs\cap\Igend$. This proves that $\Igend\neq\emptyset$.
\vskip2mm
Assume that $u_o^a=u_o^b$, with $a\neq b$, and consider the hyperplane $\Dl_{ab}:=\{u^a=u^b\}$ of $\C^n$. Denote by $\iota\colon \Dl_{ab}\cap\mathbb D(\bm u_o)\to\mathbb D(\bm u_o)$ the natural inclusion map. The germ of the integrable deformation $(\iota^*\nabla^{\rm Sab}, \iota^*\underline{\C^n},\Dl_{ab}\cap\mathbb D(\bm u_o), \bm u_o)$ defines an element of $\Ifs\setminus\Igend$. Hence, we have $\Ifs\setminus\Igend\neq \emptyset$.
\vskip2mm
It remains to show that $\Igend\subseteq \Ifs$. Consider a connection $\nabla^o$ with pole part $A_o\in\Mdiag$, and assume that Property PNR holds. Let $(\nabla,E,X,x_o)\in\Igend$ be defined by a matrix $\widehat\Om(z,x)$ as in \eqref{forcan3}, for a suitable triple of holomorphic matrices $(\Dl_0,\frak B_o',L)$ defined on an open neighborhood $U$ of $x_o$. Without loss of generality, by the Property PNR we may assume that the constant matrix $\frak B_o'={\rm diag}(b_1,\dots, b_n)$ itself is such that
\[b_i-b_j\not\in\Z\setminus\{0\}\quad \text{whenever }f_i(x_o)=f_j(x_o),\quad\text{for some }i\neq j.
\]One can always recover this condition up to a gauge equivalence $\widehat\Om\mapsto T^{-1}\widehat\Om T$ by a constant matrix $T$.
\vskip2mm
We need to show the existence of a sequence of holomorphic matrices $F_k\colon U\to M(n,\C)$, with $k\geq 1$, such that
\[F^{-1}\widehat\Om F+F^{-1}\ed F=-\ed \left(z\Dl_0(x)\right)-\frak B_o'\frac{\ed z}{z},\qquad F:={\rm Id}_n+\sum_{k=1}^\infty F_kz^{-k}. 
\]As in the proof of point (3) of Theorem \ref{mth1}, the functions $F_k$'s can be found by iteratively solving the equations
\beq
\label{8dec21.1}
[\Dl_0,F_{k+1}]+[\frak B'_o,F_k]+[L,\Dl_0]F_k+kF_k=0,\qquad k\geq 0,	\quad F_0={\rm Id}_n.
\eeq
The procedure is standard. The matrices $F_k(x)$ can be computed, entry by entry, in terms of the entries of $F_h(x)$ with $h<k$, for any $x\in U$. At each $x\in U$, the procedure works case by case, according weather $f_i(x)\neq f_j(x)$ or $f_i(x)=f_j(x)$:
\begin{itemize} 
\item if $f_i(x)\neq f_j(x)$, with $i\neq j$, then 
\beq\label{ie11}
F_{k+1}(x)_{ij}=\frac{1}{f_j(x)-f_i(x)}\left([\frak B'_o,F_k(x)]+[L(x),\Dl_0(x)]F_k(x)+kF_k(x)\right)_{ij},\quad k\geq 0;
\eeq
\item if $f_i(x)=f_j(x)$, with $i\neq j$, then 
\beq
\label{ie21}
F_{k+1}(x)_{ij}=-\frac{1}{b_i-b_j+k+1}\sum_{\ell}(f_\ell(x)-f_i(x))L(x)_{i\ell}F_{k+1}(x)_{\ell j},\qquad k\geq 0;
\eeq
\item the diagonal entries are given by
\beq\label{ie31}
F_{k+1}(x)_{ii}=-\frac{1}{k+1}\sum_{\ell\neq i}\frak B(x)_{i\ell}(F_{k+1})_{\ell i},\qquad k\geq 0.
\eeq
\end{itemize}
From this construction, it is clear that the $F_k$'s are holomorphic on $U\setminus{\rm coal}(\Dl_0)$. Moreover, on the complement $U\setminus\bigcup_{a,b}\{z\colon f_a(x)=f_b(x)\}$, the following identities hold: for any pair $(i,j)$, with $i\neq j$, we have
\beq
(\ed f_j-\ed f_i)(F_{k+1})_{ij}=\left( [F_1,\ed \Dl_0]F_{k}-\ed F_k\right)_{ij}.
\eeq
This follows from the argument used in the proof of point (3) of Theorem \ref{mth1}.
\vskip2mm
Notice that equation \eqref{8dec21.1}, specialized for $k=0$, is trivially solved by the choice $F_1''=L$. In particular, $F_1''$ holomorphically extends on the whole open set $U$, including the coalescence set ${\rm coal}(\Dl_0)$. Moreover, equation \eqref{ie31} implies that $F_1$ holomorphically extends to the whole $U$.

We claim that the holomorphicity of $F_1$ at a point $x_c\in{\rm coal}(\Dl_0)$ implies that {\it all} the matrices $F_k$ are holomorphic at $x_c$. In order to prove this, we proceed by induction on $k$. For $k=1$, the claim is tautological. Assume that $F_1,\dots, F_k$ are holomorphic on $X$. If we introduce local coordinates $(x^1,\dots, x^d)$ on $X$, with $d=\dim X$, by assumption there exists $h\in\{1,\dots, d\}$ such that $\der_hf_j-\der_h f_i$ is not vanishing on $U$. Hence, we have 
\[(F_{k+1})_{ij}=\frac{1}{\der_hf_j-\der_h f_i}\left( [F_1,\der_h\Dl_0]F_{k}-\der_hF_k\right)_{ij},
\]and the right-hand-side is holomorphic at $x_c$. This shows that the off-diagonal entries of each $F_k$'s are holomorphic at $x_c$. The diagonal entries of $F_{k+1}$ are determined by the off-diagonal ones, by equation \eqref{ie31}, and they are holomorphic at $x_c$. This completes the proof.
\endproof

\begin{rem}\label{rkest}
Due to equations \eqref{ie11} and \eqref{ie21}, the holomorphicity of $F_1$ and $F_2$ at points of ${\rm coal}(\Dl_0)$ is equivalent to the following condition on the holomorphic matrix $L$: for any $x_c\in{\rm coal}(\Dl_0)$, such that $f_i(x)=f_j(x)$ for some $i\neq j$, there exists a neighborhood $\mc U(x_c)$ where
\[
\label{holcon}(b_j-b_i-1)L_{ij}(x)-\sum_{\ell\neq i}(f_\ell(x)-f_i(x))L_{i\ell}(x)L_{\ell j}(x)=O\left(f_i(x)-f_j(x)\right),\quad x\in\mc U(x_c).
\]At this point of the presentation, it is not not obvious why such an estimate holds true. We will give a direct justification of this fact in the subsequent sections, standing on a deeper analysis of the generalized Darboux--Egoroff system of equations \eqref{DEeq1} and \eqref{DEeq2}. See Remark \ref{rkest2}.
\end{rem}

\subsection{Generalized Darboux--Egoroff equations, and its initial value property}
Consider $n$ holomorphic functions $f_1(\bm x),\dots, f_n(\bm x)$ in $d$ complex variables $\bm x=(x^1,\dots, x^d)\in\C^d$. Let $b_1,\dots, b_n\in\C$ be arbitrary constants. In what follows we set $\der_i=\frac{\der}{\der x^i}$.
\vskip2mm
The {\it generalized Darboux--Egoroff} system ${\rm DE}_{d,n}((f_i)_{i=1}^n;(b_i)_{i=1}^n)$ is the following system of PDE's for $n^2-n$ unknown functions $(F_{kh}(\bm x))_{k,h=1}^n$, with $k\neq h$:
\begin{multline}
\label{DE1}
(\DE{j}{h}{k})\der_iF_{kh}-(\DE{i}{h}{k})\der_jF_{kh}=\\
\sum_{\ell=1}^n(\DE{i}{\ell}{k})(\DE{j}{h}{\ell})F_{k\ell}F_{\ell h}-\sum_{\ell=1}^n(\DE{j}{\ell}{k})(\DE{i}{h}{\ell})F_{k\ell}F_{\ell h},
\end{multline}
\begin{multline}
\label{DE2}
(f_h-f_k)\der_iF_{kh}=(b_h-b_k-1)(\DE{i}{h}{k})F_{kh}\\+\sum_{\ell=1}^n(\DE{i}{\ell}{k})(f_h-f_\ell)F_{k\ell}F_{\ell h}
-\sum_{\ell=1}^n(f_\ell-f_k)(\DE{i}{h}{\ell})F_{k\ell}F_{\ell h},
\end{multline}
for any $i,j=1,\dots, d$, and any $k,h=1,\dots, n$, with $k\neq h$. Solutions $F$ of the generalized Darboux--Egoroff system can be arranged in a off-diagonal matrix in $M(n,\C)$. In matrix notation, equations \eqref{DE1}, \eqref{DE2} read
\[
\ed [F,\ed\Dl_0]=[F,\ed\Dl_0]\wedge [F,\ed\Dl_0],\qquad
\ed[F,\Dl_0]=\left[[F,\ed\Dl_0],\frak B'_o+[F,\Dl_0]\right],
\]where $\Dl_0(\bm x)=(f_1(\bm x),\dots, f_n(\bm x))$, and $\frak B_o'={\rm diag}(b_1,\dots, b_n)$.
\begin{rem}
In the very special case $d=n=3$, and $f_i(x^1,x^2, x^3)=x^i$, with $i=1,2,3$, the system of nonlinear PDEs \eqref{DE1} and \eqref{DE2} reduces to the generic Painlev\'e equation ${\rm PVI}_{\al,\bt,\gm,\dl}$. See \cite[Th.\,4.1]{Lor14}. 
\end{rem}
\vskip2mm
Set $\mc A:=M(n,\C)[\![(x^i-x_o^i)_{i=1}^d]\!]$.

\begin{thm}\label{mth3}
Let $\bm x_o\in\C^d$, and assume that 
\begin{itemize}
\item for any $k,h=1,\dots, n$, with $k\neq h$, we have $\ed_{\bm x_o} f_k\neq \ed_{\bm x_o} f_h$,
\item $b_h-b_k\notin\Z^*$ whenever $f_h(\bm x_o)=f_k(\bm x_o)$.
\end{itemize}
If $F_1,F_2\in \mc A$ are two formal power series solutions of the system ${\rm DE}_{d,n}((f_i)_{i=1}^n;(b_i)_{i=1}^n)$ such that $F_1(\bm x_o)=F_2(\bm x_o)$, then $F_1=F_2$. 

In particular, given an initial condition $F_o\in M(n,\C)$, with $F_o=F_o''$, there exists at most one holomorphic solution $F$ of ${\rm DE}_{d,n}((f_i)_{i=1}^n;(b_i)_{i=1}^n)$ such that $F(x_o)=F_o$.
\end{thm}

\proof
We have to show that the derivatives $\der_{i_1}\dots\der_{i_N}F_{kh}(\bm x_o)$ can be computed from the only knowledge of the numbers $F_{kh}(\bm x_o)$. We proceed by induction on $N$. Let us start with the
case $N=1$. 

Before proceeding with the proof, notice for any $k,h=1,\dots, n$, with $k\neq h$, there exists an index\footnote{Clearly $j_0$ depends on $(k,h)$, but we omit the dependence for brevity of notation.} $j_0\in\{1,\dots, d\}$ such that $\der_{j_0}f_h(\bm x_o)\neq \der_{j_0}f_k(\bm x_o)$, by assumption (1).

\vskip2mm
{\bf Step 1.} Let $a\in\{1,\dots, d\}$ be such that $\der_af_h(\bm x_o)=\der_af_k(\bm x_o)$. 
Consider the equation \eqref{DE1} with specialization of indices $(i,j)=(a,j_0)$. By evaluation at $\bm x=\bm x_o$, we can compute the number $\der_aF_{kh}(\bm x_o)$.
\vskip2mm
{\bf Step 2.} Assume that $f_k(\bm x_o)\neq f_h(\bm x_o)$. Then, for any $a\in\{1,\dots, d\}$, we can compute all the number $\der_aF_{kh}(\bm x_o)$ from equation \eqref{DE2}, with specialzation of index $i=a$, by evaluation at $\bm x=\bm x_o$.
\vskip2mm
{\bf Step 3.} Assume that $f_k(\bm x_o)= f_h(\bm x_o)$. Let $a\in\{1,\dots,d\}$ be such that $\der_af_k(\bm x_o)\neq \der_af_h(\bm x_o)$. Consider equation \eqref{DE2}, with specialization of index $i=a$, and compute the $\der_{j_0}$-derivative of both sides of the equation: we obtain
\begin{multline}
\label{intereq1}
(\DE{j_0}{h}{k})\der_aF_{kh}+(f_h-f_k)\der^2_{aj_0}F_{kh}=\\(b_h-b_k-1)(\der^2_{aj_0}f_h-\der^2_{aj_0}f_k)F_{kh}+(b_h-b_k-1)(\DE{a}{h}{k})\der_{j_0}F_{kh}\\
+\der_{j_0}\left[\sum_{\ell=1}^n(\DE{a}{\ell}{k})(f_h-f_\ell)F_{k\ell}F_{\ell h}
-\sum_{\ell=1}^n(f_\ell-f_k)(\DE{a}{h}{\ell})F_{k\ell}F_{\ell h}\right],
\end{multline}
and the last summand can be further expanded by using Leibnitz rule. By evaluating \eqref{intereq1} at $\bm x=\bm x_o$, we obtain an identity of the type
\begin{multline}
\label{lsyst1}
(\der_{j_0}f_h(\bm x_o)-\der_{j_0}f_k(\bm x_o))\der_aF_{kh}(\bm x_o)\\-(b_h-b_k-1)(\der_{a}f_h(\bm x_o)-\der_{a}f_k(\bm x_o))\der_{j_0}F_{kh}(\bm x_o)=X_1,
\end{multline}
where $X_1$ is an expression involving only 
\begin{itemize}
\item the values at $\bm x=\bm x_o$ of $f_h, f_k$, and of their first and second partial derivatives,
\item  the values at $\bm x=\bm x_o$ of $F$, and of its first partial derivatives which can be computed in Step 2.
\end{itemize}
Similarly, consider equation \eqref{DE1}, with specialization of indices $(i,j)=(a,j_0)$. By evaluation at $\bm x=\bm x_o$, we obtain an identity of the type
\beq
\label{lsyst2}
(\der_{j_0}f_h(\bm x_o)-\der_{j_0}f_k(\bm x_o))\der_aF_{kh}(\bm x_o)-(\der_{a}f_h(\bm x_o)-\der_{a}f_k(\bm x_o))\der_{j_0}F_{kh}(\bm x_o)=X_2,
\eeq
where $X_2$ is an expression involving only 
\begin{itemize}
\item the values at $\bm x=\bm x_o$ of first derivatives of $f_h, f_k$,
\item  the values at $\bm x=\bm x_o$ of $F_{kh}$.
\end{itemize}
Equations \eqref{lsyst1} and \eqref{lsyst2} define a linear system of equations 
\[\begin{pmatrix}
(\der_{j_0}f_h(\bm x_o)-\der_{j_0}f_k(\bm x_o))&-(b_h-b_k-1)(\der_{a}f_h(\bm x_o)-\der_{a}f_k(\bm x_o))\\
(\der_{j_0}f_h(\bm x_o)-\der_{j_0}f_k(\bm x_o))&-(\der_{a}f_h(\bm x_o)-\der_{a}f_k(\bm x_o))
\end{pmatrix}\begin{pmatrix}
\der_aF_{kh}(\bm x_o)\\
\der_{j_0}F_{kh}(\bm x_o)
\end{pmatrix}=\begin{pmatrix}
X_1\\
X_2
\end{pmatrix}.
\]Such a system admits a unique solution since 
\begin{multline*}\det\begin{pmatrix}
(\der_{j_0}f_h(\bm x_o)-\der_{j_0}f_k(\bm x_o))&-(b_h-b_k-1)(\der_{a}f_h(\bm x_o)-\der_{a}f_k(\bm x_o))\\
(\der_{j_0}f_h(\bm x_o)-\der_{j_0}f_k(\bm x_o))&-(\der_{a}f_h(\bm x_o)-\der_{a}f_k(\bm x_o))
\end{pmatrix}=\\
\\
=(\der_{j_0}f_h(\bm x_o)-\der_{j_0}f_k(\bm x_o))(\der_{a}f_h(\bm x_o)-\der_{a}f_k(\bm x_o))(b_h-b_k-2)\neq 0.
\end{multline*}
This proves that all the first derivatives $\der_aF_{kh}(\bm x_o)$ can be computed.
\vskip2mm
{\bf Inductive step. }Assume to know all the $N$-th derivatives $\der_{i_1}\dots\der_{i_N}F_{kh}(\bm x_o)$. We show how to compute all the $(N+1)$-th derivatives $\der_{i_1}\dots\der_{i_{N+1}}F_{kh}(\bm x_o)$.
\vskip2mm
{\bf Step 1.} Assume there exists an $\ell\in\{1,\dots, N+1\}$ such that $\der_{i_\ell}f_{h}(\bm x_o)= \der_{i_\ell}f_{k}(\bm x_o)$. Without loss of generality, we can assume $\ell=N+1$. Consider equation \eqref{DE1}, with specialization of indices $(i,j)=(i_{N+1},j_0)$, and take the $\der_{i_1}\dots\der_{i_N}$-derivative of both sides. By evaluation at $\bm x=\bm x_o$, we can compute the number $\der_{i_1}\dots\der_{i_{N+1}}F_{kh}(\bm x_o)$ in terms of lower order derivatives of $F_{kh}$ at $\bm x_o$ (hence previously computed).
\vskip2mm
{\bf Step 2.} Assume that $f_h(\bm x_o)\neq f_k(\bm x_o)$. Consider equation \eqref{DE2}, with specialization of index $i=i_{N+1}$, and take the $\der_{i_1}\dots\der_{i_N}$-derivative of both sides. By evaluation at $\bm x=\bm x_o$, we can compute the number $\der_{i_1}\dots\der_{i_{N+1}}F_{kh}(\bm x_o)$ in terms of lower order derivatives of $F_{kh}$ at $\bm x_o$ (hence previously computed).
\vskip2mm
{\bf Step 3.} Assume that $f_h(\bm x_o)= f_k(\bm x_o)$, and that for any $\ell\in\{1,\dots, N+1\}$ we have $\der_{i_\ell}f_h(\bm x_o)\neq \der_{i_\ell}f_k(\bm x_o)$. Set
\[\der_{\hat 0}:=\der_{i_1}\der_{i_2}\dots\der_{i_N+1},\quad \der_{\hat\ell}:=\der_{j_0}\der_{i_1}\dots\der_{i_{\ell-1}}\der_{i_{\ell+1}}\dots\der_{i_{N+1}},\quad \ell=1,\dots, N+1.
\]For any $\ell\in\{1,\dots, N+1\}$, consider equation \eqref{DE1}, with specialization of indices $(i,j)=(i_{\ell},j_0)$. By taking the $\der_{i_1}\dots\der_{i_{\ell-1}}\der_{i_{\ell+1}}\dots\der_{i_{N+1}}$-derivative of both sides, we obtain an identity of the form
\beq
\label{ls21}
(\DE{j_0}{h}{k})\der_{\hat 0}F_{kh}-(\DE{i_{\ell}}{h}{k})\der_{\hat{\ell}}F_{kh}=Z_{1,\ell},
\eeq
where $Z_{1,\ell}$ is a polynomial expression in the $p$-th derivatives of $(f_h-f_k)$, with $0\leq p\leq N+1$, and the $q$-th derivatives of $F$ with $0\leq q\leq N$.

Consider equation \eqref{DE2}, with specialization of index $i=i_{N+1}$. By taking the $\der_{\widehat{N+1}}$-derivative of both sides we obtain an identity of the form
\begin{multline}
\label{ls22}
(f_h-f_k)\der^{N+2}_{j_0i_1\dots i_{N+1}}F_{kh}+(\DE{j_0}{h}{k})\der_{\hat 0}F_{kh}+\sum_{\ell=1}^{N}(\DE{i_\ell}{h}{k})\der_{\hat\ell}F_{kh}\\-(b_h-b_k-1)(\DE{i_{N+1}}{h}{k})\der_{\widehat{N+1}}F_{kh}=Z_2,
\end{multline}where $Z_2$ is a polynomial expression in 
\begin{itemize}
\item the $p$-th derivatives of $f_h-f_k$, with $0\leq p\leq N+2$, 
\item in the $q$-th derivatives of $F$, with $0\leq q\leq N+1$.
\end{itemize}Moreover, by evaluating at $\bm x=\bm x_o$ both sides of \eqref{ls22}, one can notice that:
\begin{itemize}
\item the first term in the left-hand-side of \eqref{ls22} (i.e.\,\,the one with the $(N+2)$-th derivative of $F_{kh}$) cancels,
\item the only $(N+1)$-th derivatives of $F$ appearing in $Z_2(\bm x_o)$ are those computed in Step 2. 
\end{itemize}
Hence, equations \eqref{ls21} evaluated at $\bm x=\bm x_o$, for $\ell=1,\dots, N+1$, and equations \eqref{ls22} evaluated at $\bm x=\bm x_o$, define a linear system of equations in the $N+2$ unknowns $\der_{\hat 0}F_{kh}(\bm x_o),$ $\der_{\hat 1}F_{kh}(\bm x_o),$ $\dots, \der_{\widehat{N+1}}F_{kh}(\bm x_o)$:
\[W\begin{pmatrix}
\der_{\hat 0}F_{kh}(\bm x_o)\\
\der_{\hat 1}F_{kh}(\bm x_o)\\
\vdots\\
\der_{\widehat{N}}F_{kh}(\bm x_o)\\
\der_{\widehat{N+1}}F_{kh}(\bm x_o)
\end{pmatrix}=
\begin{pmatrix}
Z_{1,1}(\bm x_o)\\
Z_{1,2}(\bm x_o)\\
\vdots\\
Z_{1,N+1}(\bm x_o)\\
Z_2(\bm x_o)
\end{pmatrix},
\]where the matrix $W$ equals
\[W=\begin{pmatrix}
D_0&-D_1&0&0&\dots&0&0\\
D_0&0&-D_2&0&\dots&0&0\\
D_0&0&0&-D_3&\dots&0&0\\
\vdots&\vdots&\vdots&\vdots&\ddots&\vdots&\vdots\\
D_0&0&0&0&\dots&-D_N&0\\
D_0&0&0&0&\dots&0&-D_{N+1}\\
D_0&D_1&D_2&D_3&\dots&D_N&-\kappa D_{N+1}
\end{pmatrix},\qquad 
\begin{array}{c}
D_0:=\der_{j_0}f_h(\bm x_o)-\der_{j_0}f_k(\bm x_o),\\
\\
D_\ell:=\der_{i_\ell}f_h(\bm x_o)-\der_{i_\ell}f_k(\bm x_o),\\
\\
\kappa=b_h-b_k-1.
\end{array}
\]We have
\[\det W=(\der_{j_0}f_h(\bm x_o)-\der_{j_0}f_k(\bm x_o))\prod_{\ell=1}^{N+1}(\der_{i_\ell}f_h(\bm x_o)-\der_{i_\ell}f_k(\bm x_o))(N+1-b_h+b_k)\neq 0.
\]This proves that all the $(N+1)$-th derivatives $\der_{i_1}\dots\der_{i_{N+1}}F_{kh}(\bm x_o)$ can be computed.
\endproof

\subsection{$\frak I$-universal integrable deformations}Consider a meromorphic connection $\nabla^o$ on a vector bundle $E^o\to \Pb^1$, with matrix of connection 1-forms 
\[\Om_o=-\left(A_o+\frac{1}{z}B_o\right)\ed z.
\]Assume that one of the following assumptions holds true:
\begin{enumerate}
\item[(I)] $A_o\in\Mdiag\cap\Mreg$;
\item[(II)] $A_o\in\Mdiag$, and the Property PNR is satisfied.
\end{enumerate}

\begin{defn}\label{iun}
Let $\frak I$ be a class of integrable deformations of $\nabla^o$. An integrable deformation $(\nabla,E,X,x_o)$ of $\nabla^o$ is ${\frak I}$-{\it versal}, if 
\begin{itemize}
\item $(\nabla,E,X,x_o)$ is an element of $\frak I$,
\item any element $(\nabla',E',X',x'_o)$ of $\frak I$ is induced by $(\nabla,E,X,x_o)$ via pull-back along a base-change $\phi\colon (X',x_o')\to (X,x_o)$.
\end{itemize}
It is $\frak I$-{\it universal}, if the germ at $x_o'$ of the base change $\phi$ is unique. 
\end{defn}

\begin{rem}
Given a integrable deformation $\nabla$ of the connection $\nabla^o$, there always exist classes $\frak I$ such that $\nabla$ is $\frak I$-versal. For example, $\nabla$ is clearly $\{\nabla\}$-versal. This is not true for universality, due to possible ``internal symmetries'' of the integrable deformations. Here is an example. Let $f\colon\C\to\C$ be an arbitrary holomorphic function, and $d_1,d_2\in\C$. Consider the matrices $\Dl_0,\frak B\colon\C^2\to M(2,\C)$ and the 1-form $\om\colon\C^2\to M(2,\C)\otimes \Om^1_\C$ defined by
\[\Dl_0(x)=\begin{pmatrix}
f(x_1^2+x_2^2)&0\\
0&f(x_1^2+x_2^2)
\end{pmatrix},\quad 
\frak B(x)=\begin{pmatrix}
d_1&(d_1-d_2)(x_1^2+x_2^2)\\
0&d_2
\end{pmatrix},
\]
\[\om(x)=\begin{pmatrix}
0&2\\
0&0
\end{pmatrix}(x_1\ed x_1+x_2\ed x_2).
\]The matrices 
\[\widetilde\Om_o(z)=-\left(\Dl_0(0)+\frac{1}{z}\frak B(0)\right)\ed z, \quad \widetilde \Om(z,x)=-\left(\Dl_0(x)+\frac{1}{z}\frak B(x)\right)\ed z-z\,\ed\Dl_0(x)+\om,
\]define two connections $\nabla^o$ and $\nabla$, respectively. The connection $\nabla$ is an integrable deformation of $\nabla^o$, at which it restricts at $x=0$. The connection $\nabla$ is not $\{\nabla\}$-universal: any linear map $l\colon\C^2\to\C^2$, defined by a matrix in $O(2,\C)$, is such that $l^*\widetilde\Om=\widetilde\Om$.
\end{rem}

Let us denote\footnote{Here the superscript ``JMUMS'' stands for Jimbo--Miwa--Ueno--Malgrange--Sabbah.} by $(\nabla^{\rm JMUMS}, \underline{\C^n}, \mathbb D,\bm u_o)$ the integrable deformation of $\nabla^o$ whose existence is guaranteed by Theorems \ref{thmal2} and \ref{thmal} (Case I), and Theorem \ref{tsab} (Case II). Here $\mathbb D\subseteq\C^n$ denotes a sufficiently small polydisc centered at $\bm u_o$. Recall that $\nabla^{\rm JMUMS}$ has matrix of connection 1-forms given by
\[\widehat\Om_{\rm JMUMS}(z,\bm u)=-\ed\left(z\La(\bm u)\right)-\left([\Gm(\bm u),\La(\bm u)]+\mc B'_o\right)\frac{\ed z}{z}-[\Gm(\bm u),\ed \La(\bm u)],
\]where $\La\colon \C^n\to M(n,\C)$, is defined by $\La(\bm u)={\rm diag}(u^1,\dots, u^n)$, and where the matrix $\mc B'_o={\rm diag}(b_1,\dots, b_n)$ is constant, and satisfying Property PNR-$(\star\star)$ in Case II. 

\begin{thm}\label{mth4}Let $\frak I$ be a class of integrable deformations of the connection $\nabla^o$ satisfying conditions (I) or (II) above.
\begin{enumerate}
\item An integrable deformation of $\nabla^o$ is induced by $(\nabla^{\rm JMUMS}, \underline{\C^n}, \mathbb D,\bm u_o)$ only if it is of fs-type.
\item If the deformation $(\nabla^{\rm JMUMS}, \underline{\C^n}, \mathbb D,\bm u_o)$ is $\frak I$-versal, then it is $\frak I$-universal.
\item There exists a unique maximal class $\frak I_{\rm JMUMS}$ of integrable deformations of $\nabla^o$ such that $\nabla^{\rm JMUMS}$ is $\frak I_{\rm JMUMS}$-universal.
\item In Case I, we have $\frak I_{\rm JMUMS}=\frak I(\nabla^o)$.
\item In Case II, we have $\Igend\subseteq{\frak I_{\rm JMUMS}}\subseteq \Ifs$.
\end{enumerate}
\end{thm}

\proof Point (1) follows from point (1) of Theorem \ref{mth1}, since $\nabla^{\rm JMUMS}\in\Ifs$.
\vskip2mm
Assume that $(\nabla,E,X,x_o)\in\Ifs$ is induced by $(\nabla^{\rm JMUMS}, \underline{\C^n}, \mathbb D,\bm u_o)$ via a holomorphic map $\phi\colon (X,x_o)\to (\C^n,\bm u_o)$, $x\mapsto (\phi_1(x),\dots, \phi_n(x))$. On the one hand, by Theorems \ref{thforcan1}, \ref{thDL}, and \ref{mth1}, the pole part of $\nabla$ is locally holomorphically diagonalizable at $x_o$, with holomorphic diagonal form $\Dl_0(x)=(f_1(x),\dots, f_n(x))$. On the other hand, the connection $\nabla$ has then matrix of connection 1-forms $\phi^*\widehat\Om_{\rm JMUMS}$: in particular, its pole part has Jordan diagonal form
\[(\phi^*\La)(x)={\rm diag}(\phi^1(x),\dots, \phi^n(x)),\quad x\in X.
\]Hence, the germ of the map $\phi$ at $x_o$ is uniquely determined: it is given by the {\it spectrum map} $\si\colon X\to \C^n$, $x\mapsto (f_1(x),\dots, f_n(x))$. This proves point (2).
\vskip2mm
The class $\frak A=\{\frak I\subseteq \frak I(\nabla^o)\colon \nabla^{\rm JMUMS}\text{ is $\frak I$-universal}\}$ is non-empty, since $\{\nabla^{\rm JMUMS}\}\in\frak A$. Consider the poset $(\frak A,\subseteq)$. By taking unions, it is easy to see that 
\begin{enumerate}
\item[(i)] every chain in $\frak A$ has a maximal element,
\item[(ii)] and $\frak A$ is upward-directed (i.e.\,\,given $\frak I_1,\frak I_2\in\frak A$, there exists $\frak I_3\in\frak A$ such that $\frak I_1,\frak I_2\subseteq\frak I_3$).
\end{enumerate}
By Zorn Lemma, $\frak A$ has a maximal element. It necessarily is unique, by (ii). This proves point (3).
\vskip2mm
In Case I, the statement $\frak I_{\rm JMUMS}=\frak I(\nabla^o)$ is equivalent to the universality predicated in Theorems \ref{thmal2} and \ref{thmal}. Hence (4) holds.
\vskip2mm
The only non-tivial inculsion of point (5) is $\Igend\subseteq\frak I^{\rm JMUMS}$. Consider an integrable deformation $(\nabla,E,X,x_o)\in\Igend$. By Theorem \ref{mth1}, $\nabla$ is of dv-type, and it has a matrix of connection 1-forms 
\[\widehat\Om(z,x)=-\left(\Dl_0(x)+\frac{1}{z}(\mc B'_o+[L(x),\Dl_0])\right)\ed z-z\ed\Dl_0(x)-[L(x),\ed\Dl_0(x)],
\]where $\mc B'_o$ satisfies Property PNR-$(\star\star)$, and $(\Dl_0,L)$ is a suitable pair of holomorphic matrices defined in a neighborhood of $x_o$. Consider spectrum map $\si\colon (X,x_o)\to (\C^n,\bm u_o)$ defined by $x\mapsto (f_1(x),\dots, f_n(x))$. Since $\si^*\La=\Dl_0$, we have
\[\widehat\Om-\si^*\widehat\Om_{\rm JMUMS}=-[L-\si^*\Gm,\Dl_0]\frac{\ed z}{z}-[L-\si^*\Gm,\ed \Dl_0].
\]Moreover, since the restriction at $x_o$ of both $\widehat\Om$ and $\si^*\widehat\Om_{\rm JMUMS}$ equals \eqref{con1.2}, we deduce that
\[L(x_o)-(\si^*\Gm)(x_o)\in \mc K\left[[-,\Dl_0];x_o\right]
\equiv\ker[-,\ed\Dl_0(x_o)],\qquad\text{that is }L(x_o)=(\si^*\Gm)(x_o).
\]Since both $L$ and $\si^*\Gm$ solves the generalized Darboux--Egoroff equations \eqref{DEeq1}, \eqref{DEeq2}, we conclude that $L=\si^*\Gm$, by Theorem \ref{mth3}. This shows that $\nabla=\si^*\nabla^{\rm JMUMS}$. This completes the proof.
\endproof

\begin{rem}\label{rkest2}
From the equality $L=\si^*\Gm$, we are able to justify the estimate of Remark \ref{rkest}. Consider equation \eqref{DE03}: by evaluating both sides at $\bm u(x)=(f_1(x),\dots, f_n(x))$, we obtain 
\begin{multline}
\label{holcon}(b_j-b_i-1)L_{ij}(x)-\sum_{\ell\neq i}(f_\ell(x)-f_i(x))L_{i\ell}(x)L_{\ell j}(x)=\left(f_i(x)-f_j(x)\right)\cdot \der_j\Gm_{ij}|_{\bm u(x)}\\
=O\left(f_i(x)-f_j(x)\right),
\end{multline}
as expected from the proof of Theorem \ref{mth2}.
\end{rem}

The following result provides a further description of the class $\frak I_{\rm JMUMS}$.
\begin{prop} Let $\nabla^o$ satisfy condition (I) or (II). Consider an integrable deformation $(\nabla,E,X,x_o)$ in $\Ifs$, with pole part $\Dl_0(x)=(f_1(x),\dots, f_n(x))$. 
\begin{enumerate}
\item 
The connection $\nabla$ is formally gauge equivalent to the connection $\si^*\nabla^{\rm JMUMS}$, where $\si\colon X\to \C^n$ is the spectrum map $x\mapsto (f_1(x),\dots, f_n(x))$. This means that there exist holomorphic functions $\Phi_k\colon X\to M(n,\C)$, with $\Phi_0(X)\subseteq GL(n,\C)$, such that
\[
\Phi^{-1}\widehat\Om\Phi+\Phi^{-1}\ed \Phi=\si^*\widehat\Om_{\rm JMUMS},\qquad \Phi(z,x)=\sum_{k\geq 0}\Phi_k(x)z^{-k}.
\] Moreover, if $\Phi_0(x)={\rm Id}_n$ for any $x\in X$, such a formal gauge equivalence is unique.
\item We have $\nabla\in \frak I_{\rm JMUMS}$ if and only if the formal gauge equivalence above is actually convergent.
\end{enumerate}
\end{prop}
\proof
Both $\nabla$ and $\si^*\nabla^{\rm JMUMS}$ are formally simplifiable: they are formally equivalent to the connection $\ed-\ed \left(z\Dl_0(x)\right)-\frak B_o'\frac{\ed z}{z}$, via unique formal gauge equivalences of the form ${\rm Id}_n+O(z^{-1})$. Point (1) follows.

If $\nabla\in \frak I_{\rm JMUMS}$, then there exists an analytic gauge equivalence of the form $T(z,x)={\rm Id}_n+\sum_{k\geq 1}T_k(x)z^{-k}$ such that $T^{-1}\widehat\Om T+T^{-1}\ed T=\si^*\widehat\Om_{\rm JMUMS}$. By uniqueness, we have $T=\Phi$. Conversely, if $\Phi$ is convergent, then $\nabla$ and $\si^*\nabla^{\rm JMUMS}$ are analytically gauge equivalent.
\endproof

\appendix
\section{}

\subsection{Proof of Theorem \ref{thmal1}}\label{appmal} Consider a trivial vector bundle $E^o$ on $\Pb^1$, equipped with a meromoprhic connection $\nabla^o$, admitting (in a suitable basis of sections) a matrix of 1-forms connection of the form \eqref{con1} with $A_o$ regular.  Let $(E,\nabla,\bm x_o)$ be an arbitrary integrable deformation of $(E^o,\nabla^o)$, parametrized by a simply connected manifold $X$. By Theorem \ref{thimp}, the deformation admits a matrix of connection 1-forms as in \eqref{condef}, whose coefficients are subjected to the integrability equations \eqref{inteq1}. 
\vskip2mm
Since $A_o$ is regular, the matrix $A(\bm x)$ is regular for $\bm x$ is a sufficiently small neighborhood $U\subseteq X$ of $\bm x_o$. From the equation $[A,C]=0$, we deduce that $C(\bm x)$ is a polynomial expression of $A(\bm x)$ for $\bm x\in U$. Hence, the equations $C\wedge C=0$ are automatically satisfied for $\bm x\in U$.

Since $\ed C=0$, at least locally we have $C=\ed K$ for some holomorphic matrix-valued function $K$: up to adding a constant matrix, we can assume that $K(x_o)=0$. Consequently, system \eqref{inteq1} is equivalent to
\[A-K-[K,B_o]={\rm const.}=A_o,\qquad \qquad [A,\ed K]=0.
\]These equations define a Pfaffian system for the matrix $K$ only:
\beq\label{pfaff}
\om=[A_o+K+[K,B_o],\ed K]=0.
\eeq
Given a solution $K$, the matrix $A$ can be reconstructed by $A=A_o+K+[K,B_o]$.
\vskip2mm
The argument above shows that there is a 1-1 correspondence between integrable deformations of a connection $(E^o,\nabla^o)$ and germs of maps $\varphi\colon (X,x_o)\to (M(n,\C),0)$ such that $\varphi^*\om=0$. The maps $\phi$ are the integral manifolds of the Pfaffian system \eqref{pfaff} passing through $0\in M(n,\C)$.

\begin{lem}
If $A_o$ is regular, then the Pfaffian system \eqref{pfaff} is completely integrable on $M(n,\C)$. The maximal integral solutions define a foliation: in a neighborhood of $0\in M(n,\C)$, the leaves are $n$-dimensional. 
\end{lem}

\proof
Let $v_1,v_2$ be two vector fields on $M(n,\C)$ defined in a neighborhood of the origin. We have to show that if $\om(v_i)=0$ for $i=1,2$, then also $\ed\om(v_1\wedge v_2)=0$. A simple computation shows that $\ed\om=2\ed K\wedge \ed K+2[\ed K\wedge \ed K, B_o]$. By identifying $v_1,v_2$ with matrices in $M(n,\C)$, the condition $\om(v_i)=0$ is equivalent to $[A(x),v_i]=0$, for $i=1,2$. For $x\in U$ as above, the matrix $A(x)$ is regular, and the matrices $v_1,v_2$ are polynomials in $A(x)$. In particular, we have $[v_1,v_2]=0$, so that $2[v_1,v_2]+2[[v_1,v_2],B_o]=0.$ This is exactly the condition $\ed\om(v_1,v_2)=0$. Finally, notice that the dimension of the leaf passing through $0$ equals the dimension of the centralizer of $A_o$ in $M(n,\C)$. This equals $n$, since $A_o$ is regular.
\endproof

The germ of the universal deformation of $(E^o,\nabla^o)$ is given by the germ of the maximal integral submanifold of the Pfaffian system \eqref{pfaff} passing through $0$. This completes the proof of Theorem \ref{thmal1}.

\subsection{Versal deformations do not exist if $A_o\notin \Mreg$}\label{appnover} Let $n=2$, and introduce the matrices 
\beq\label{eqap1}A_o=\begin{pmatrix}
0&0\\
0&0
\end{pmatrix},\quad B_o=\begin{pmatrix}
c&0\\
0&0
\end{pmatrix}.
\eeq Let $\nabla^o$ be the connection on $\underline{\C^2}\to \Pb^1$ with connection \eqref{con1}. Consider an integrable deformation of $\nabla^o$ defined by a matrix  $\Om(z,\bm x)$ as in \eqref{condef}, where $\bm x$ is a parameter varying in a polydisc $\mathbb D\subseteq \C^m$. The matrix
\[\Om'(z,\bm x,s):=\Om(z,\bm x)+\ed(zs\,{\rm Id}_2)
\]defines a new integrable deformation of $\nabla^o$, parametrized by points $(\bm x,s)\in\mathbb D\times \C$. If $A(\bm x)$ is the pole part of $\Om$, consider the hypersurface $\mc L:=\{s=\frac{1}{2}{\rm Tr}\,A(\bm x)\}$ in $\mathbb D\times \C$: the restriction $\Om'|_{\Pb^1\times\mc L}$ is a deformation of $\nabla^o$ with traceless pole part. This shows that, without loss of generality, we may restrict to the study of integrable deformations with traceless pole part \[A(\bm x)=\begin{pmatrix}
\al(\bm x)&\bt(\bm x)\\
\gm(\bm x)&-\al(\bm x)
\end{pmatrix}.\]

As in the previous section, the integrability system \eqref{inteq1} can be reduced to the following system of equations in the pair $(A,K)$, with $\ed K=C$:
\beq
\label{pf2}
[A,\ed K]=0,\quad A=K+[K,B_o],\quad \ed K\wedge \ed K=0.
\eeq
Assume that $c\neq \pm 1$, so that the operator \[\nu\colon M(2,\C)\to M(2,\C),\quad X\mapsto X+[X,B_o],\] is invertible, with inverse 
\[\nu^{-1}\colon \begin{pmatrix}
x_{11}&x_{12}\\
x_{21}&x_{22}
\end{pmatrix}\mapsto \begin{pmatrix}
x_{11}&\frac{1}{1-c}x_{12}\\
\frac{1}{1+c}x_{21}&x_{22}
\end{pmatrix}.
\]The system \eqref{pf2} can be reduced to a system of differential equations in $A$ only, since $K=\nu^{-1}(A)$. The equation $[A,\ed K]=0$ becomes 
\[\begin{pmatrix}
0&-2\bt\\
2\gm&0
\end{pmatrix}\ed\al +
\begin{pmatrix}
-\gm&2\al\\
0&\gm
\end{pmatrix}\frac{\ed\bt}{1-c}+
\begin{pmatrix}
\bt&0\\
-2\al&-\bt
\end{pmatrix}\frac{\ed\gm}{1+c}=0,
\]which defines following the Pfaffian system on the space $\C^3$ of triples $(\al,\bt,\gm)$: 
\begin{align}\label{pf2.2}
\nonumber
\om_1&=(1-c)\bt\,\ed\al-\al\,\ed\bt=0,\\
\om_2&=(1+c)\gm\,\ed\al-\al\,\ed\gm=0,\\
\nonumber
\om_3&=(1+c)\gm\,\ed\bt-(1-c)\bt\,\ed\gm=0.
\end{align}This system can be written in a more compact way as
\[
\frac{\ed\al}{\al}=\frac{\ed\bt}{(1-c)\bt}=\frac{\ed\gm}{(1+c)\gm}.
\]
The equation $\ed K\wedge \ed K=0$ reduces to
\[\ed\al\wedge \ed\bt=\ed\al\wedge \ed\gm=\ed\bt\wedge \ed\gm=0,
\]which are automatically satisfied if \eqref{pf2.2} holds true.

The discussion in the previous section shows that integrable deformations of $\nabla^o$, with traceless pole part, define (and are defined by) germs of maps $\phi\colon (X,\bm x_o)\to (\C^3,A_o)$ such that $\phi^*\om_i=0$ for $i=1,2,3$, i.e.\,\,integral submanifolds of the Pfaffian system \eqref{pf2.2} passing through $(\al_o,\bt_o,\gm_o)=0$. 
\vskip2mm
If $(\al_o,\bt_o,\gm_o)\neq 0$, there exists a unique (one dimensional) maximal integral submanifold of the Pfaffian system above passing through $(\al_o,\bt_o,\gm_o)$:
\[\al(t)=\al_oe^t,\qquad \bt(t)=\bt_oe^{(1-c)t},\qquad \gm(t)=\gm_oe^{(1+c)t}.
\]

If $(\al_o,\bt_o,\gm_o)=0$, on the other hand, one can find many integral curves passing through $(\al_o,\bt_o,\gm_o)$. For example:
\begin{align*}
\al=t,\quad \bt=\gm=0;\qquad
\al=\gm=0,\quad \bt=t;\qquad
\al=\bt=0,\quad \gm=t;
\end{align*}
these lines are not contained in a surface, hence there is no versal deformation inducing all of them.
\vskip2mm
Infinite families of solutions arise if $c\in\mathbb Q$. Assume $c=\frac{p}{q}$ with $(p,q)=1$, and $q>0$:
\begin{itemize}
\item if $c<-1$ (i.e.\,\,$p<-q$) we have the infinite family 
\[\al=\al_0t^q,\qquad \bt=\bt_0t^{q-p},\qquad \gm=0,\quad \al_0,\bt_0\in\C;
\]
\item if $-1<c<1$ (i.e.\,\,$-q<p<q$) we also have another infinite family
\[\al=\al_0t^q,\qquad \bt=\bt_0t^{q-p},\qquad \gm=\gm_0t^{p+q},\quad \al_0,\bt_0,\gm_0\in\C;
\]if moreover $p+q$ is even, then we also have the infinite family
\[\al=0,\qquad \bt=\bt_0t^{\frac{q-p}{2}},\qquad \gm=\gm_0t^\frac{p+q}{2},\quad \bt_0,\gm_0\in\C;
\] 
\item if $c>1$ (i.e.\,\,$p>q$) we have the infinite family
\[\al=\al_0t^q,\qquad \bt=0,\qquad \gm=\gm_0t^{p+q},\quad \al_0,\gm_0\in\C.
\]
\end{itemize}
For a complete classification of the solutions see the recent paper \cite[Sec.\,8]{Her21}.

\subsection{Case of $\mc B_o'$ partially resonant}\label{nPNR} Let us now consider the case of the matrices $A_o,B_o$ as in \eqref{eqap1} with $c=1$. For such a pair of matrices the Property PNR cannot hold true. As before, let $\nabla^o$ be the connection on $\underline{\C^2}\to \Pb^1$ with connection \eqref{con1}, and consider an arbitrary (germ of) integrable deformation $\nabla$ of $\nabla^o$, parametrized by (the germ of) a pointed manifold $(X,\bm x_o)$. Let the matrix $\Om$ of 1-forms of $\nabla$ to be 
\[\Om(z,\bm x)=-\left(A(\bm x)+\frac{1}{z}B_o\right)\ed z-zC(\bm x),
\]as in Theorem \ref{thimp}.
\begin{prop}
The matrix $\Om$ is of one, and only one, of the following types:\newline
\noindent {\bf Type I:} There exist two holomorphic functions $g,m\colon X\to \C$, not identically equal, with $g(\bm x_o)=m(\bm x_o)=0$, and a complex \emph{modulus} $\kappa\in\C$ such that
\beq\label{tI}
A=\begin{pmatrix}
g&0\\
2\kappa(m-g)^2&m
\end{pmatrix},\qquad C=\begin{pmatrix}
\ed g&0\\
\kappa\cdot\ed (m-g)^2&\ed m
\end{pmatrix}.
\eeq
\noindent{\bf Type II:} There exist a holomorphic function $g\colon X\to \C$, with $g(\bm x_o)=0$, such that
\beq\label{tIb}
A=\begin{pmatrix}
g&0\\
0&g
\end{pmatrix},\qquad C=\begin{pmatrix}
\ed g&0\\
0&\ed g
\end{pmatrix}.
\eeq
\noindent{\bf Type III:} There exist two holomorphic functions $g,h\colon X\to \C$, with $h$ not identically zero, and $g(\bm x_o)=h(\bm x_o)=0$, such that
\beq\label{tII} A(\bm x)=
\begin{pmatrix}
g&0\\
0&g
\end{pmatrix},\qquad C(\bm x)=\begin{pmatrix}
\ed g&\ed h\\
0&\ed g
\end{pmatrix}.
\eeq
\end{prop}
\proof
The integrability conditions \eqref{inteq1} can be put in the form \eqref{pf2}, but we cannot reduce anymore the problem to a Pfaffian system for $A(\bm x)$, as we did in the previous section, since $c=1$. So, let 
\[K(\bm x)=\begin{pmatrix}
g(\bm x)&h(\bm x)\\
\ell(\bm x)&m(\bm x)
\end{pmatrix},
\]for arbitrary functions $g,h,\ell,m\colon X\to\C$ all vanishing at $\bm x_o$, and let us consider the system \eqref{pf2} in the variable $K$. The first equation of the system \eqref{pf2} is then equivalent to the set of equations
\beq\label{defeqdeg}
\ell\,\ed h=0,\qquad (g-m)\ed h=0,\qquad (m-g)\ed \ell-2\ell\,\ed(m-g)=0.
\eeq
We may have two cases: (Case I) either $h$ is identically zero on $X$, (Case II) or $h(\bm x)\neq 0$ for at least one $\bm x\in X$. 
\vskip2mm
In Case I, from \eqref{defeqdeg} we deduce that $\ell=\kappa (m-g)^2$ for some $\kappa\in\C$, and the remaining equation $\ed K\wedge\ed K=0$ is automatically satisfied. Consequently, the matrix of connection 1-forms is necessarily of the form \eqref{tI} or \eqref{tIb}. 
In case II, from \eqref{defeqdeg} we deduce that $\ell=0$ and $g=m$, and the remaining equation $\ed K\wedge\ed K=0$ is automatically satisfied.
Hence, the matrix $\Om$ of connection 1-forms is of the form \eqref{tII}. \endproof

\begin{prop}\label{propA3}$\quad$
\begin{enumerate}
\item All integrable deformations in $\Id$ are of Type I or II. Moreover, any element of $\Igend$ is of Type I.
\item Integrable deformations of Type I cannot induce deformations of Type III, and vice-versa. 
\item Both deformations of Type I and III can induce deformations of Type II.
\item 
Two deformations $(\nabla^{[1]},\kappa_1),(\nabla^{[2]},\kappa_2)$ of Type I cannot be induced by a same deformation if $\kappa_1\neq \kappa_2$.
\item Integrable deformations of Type I and fixed modulus $\kappa$ are induced by a universal deformation of such a type: it is the connection on $\underline{\C^2}\to\Pb^1\times \C^2$ with matrix of connection 1-forms
\[\Om^{[{\rm I}]}_\kappa(z,\bm u)=-\ed\left(z\La(\bm u)\right)-\left([F_\kappa,\La(\bm u)]+B_o\right)\frac{\ed z}{z}-[F_\kappa,\ed \La(\bm u)],
\]where
\[\La(\bm u)={\rm diag}(u_1,u_2),\quad F_\kappa=\begin{pmatrix}
0&0\\
-2\kappa&0
\end{pmatrix}.
\]
\item Integrable deformations of Type II do not admit a versal deformation. 
\item Integrable deformations of Type III admit a versal (but not universal) deformation of such a type: it is the connection on $\underline{\C^2}\to\Pb^1\times \C^2$ with matrix of connection 1-forms
\[\Om^{[{\rm III}]}(z,\bm u)=-\left(\begin{pmatrix}
u_1&0\\
0&u_1
\end{pmatrix}+\frac{1}{z}B_o\right)\ed z-z\begin{pmatrix}
\ed u_1&\ed u_2\\
0&\ed u_1
\end{pmatrix}.
\]
\end{enumerate}
\end{prop}

\proof In all Types I, II, and III,  the matrix $A(\bm x)$ is holomorphically diagonalizable: for Type I, indeed, we have
\[M^{-1}AM={\rm diag}(g(\bm x),m(\bm x)),\qquad M(\bm x)=\begin{pmatrix}
1&0\\
2\kappa (g(\bm x)-m(\bm x))&1
\end{pmatrix};
\]in Type II and III, the matrix $A$ is already in diagonal form. However, the deformation part $C$ is diagonalizable only in Types I and II. In Type I we have $M^{-1}CM={\rm diag}(\ed g,\ed m)$. Moreover, any element of $\Igend$ necessarily is of Type I. This proves point (1).

All the remaining statements easily follow from the explicit equations \eqref{tI}, \eqref{tIb}, and \eqref{tII}. 
\endproof

Points (1) and (4) of Proposition \ref{propA3} imply that there exist no $\Igend$-versal integrable deformations.

\end{document}